\documentclass[%
 aip,
 amsmath,amssymb,
preprint,
10pt]{revtex4-1}

\usepackage{graphicx}
\usepackage{dcolumn}
\usepackage{bm}

\usepackage[utf8]{inputenc}
\usepackage[T1]{fontenc}
\usepackage{mathptmx}

 \date{Revised version of \today}

\usepackage{amssymb}

\usepackage{amscd}
\usepackage{amsthm}
\usepackage[all]{xy}
\usepackage{extarrows}
\usepackage{rotating}
\usepackage[colorlinks,final,backref=page,hyperindex]{hyperref}

\newtheorem{Thm}{Theorem}[section]
\newtheorem{Lem}[Thm]{Lemma}
\newtheorem{Def}[Thm]{Definition}

\newtheorem{Prop}[Thm]{Proposition}
\newtheorem{Prop-Def}[Thm]{Proposition-Definition}

\newtheorem{Rem}[Thm]{Remark}

\newenvironment{Proof}[1][Proof]{\begin{trivlist}
\item[\hskip \labelsep {\bfseries #1}]}{\flushright
$\Box$\end{trivlist}}

\renewcommand{\k}{\Bbbk}

\renewcommand{\k}{\mathbb{K}}
\newcommand{\PH}{\operatorname{HP}}
\newcommand{\rmd}{\mathrm{d}}

\newcommand{\lra}{\longrightarrow}

\newcommand{\ra}{\rightarrow}
\newcommand{\sdp}{\times\kern-.2em\vrule height1.1ex depth-.05ex}
\newcommand{\epi}{\lra \kern-.8em\ra}

\newcommand{\px}{\frac{\partial}{\partial x}}
\newcommand{\py}{\frac{\partial}{\partial y}}

\newcommand{\sn}{{\mathrm{\mathbf{SN}}}}

\begin{document}

\preprint{AIP/123-QED}

\title[Poisson cohomology of plane Poisson structures with isolated singularities  revisited]{Poisson cohomology of plane Poisson structures with isolated singularities  revisited}

\author{Zihao Qi}
 \email{qizihao@foxmail.com}
 \affiliation{
School of Mathematical Sciences,
Shanghai Key laboratory of PMMP,
East China Normal University,\\
Shanghai 200241, People's Republic of China}
\author{Guodong Zhou}
\email[Author to whom correspondence should be addressed:]{gdzhou@math.ecnu.edu.cn}
\affiliation{
School of Mathematical Sciences,
Shanghai Key laboratory of PMMP,
East China Normal University,\\
Shanghai 200241, People's Republic of China}

\date{Revised version of \today}

\keywords{GAGA type results, Gerstenhaber algebra structure,  isolated singularities, Poisson cohomology, simple singularities\\
Mathematics Subject Classification(2020): 17B63, 14H20, 53D17.}

\begin{abstract}  Continuing   a work of   Ph.~Monnier, we determine  the Gerstenhaber algebra structure over the  Poisson cohomology groups for a large  class of   Poisson structures with isolated singularities over the plane. It reveals that there exists a  GAGA type phenomenon. We give an explicit  description via generators and relations for the case of simple singularities.

\end{abstract}

\maketitle

\section{Introduction}

The Poisson cohomology   of Poisson manifolds (and of Poisson algebras) was introduced by A.~Lichnerowicz
 \cite{Lichnerowicz}.
 J.~Huebschmann
 studied this cohomology theory from the viewpoint of Lie-Rinehart algebras \cite{Hue90, Hue98, Hue99}.
 While Poisson cohomology is  ubiquitous in Poisson geometry,  their computation is  in general very difficult.
  It has been  shown      that for symplectic manifolds, Poisson cohomology is naturally isomorphic to de Rham cohomology \cite{Lichnerowicz}.   Cohomology of regular Poisson manifolds was investigated by P.~Xu \cite{Xu92},  and  V.~Ginzburg and A.~Weinstein studied cohomology of Poisson-Lie groups \cite{GW}.

  In this paper, we are interested into cohomology of   Poisson structures over the plane.
   C.~Roger and   P.~Vanhaecke \cite{RP} computed the Poisson cohomology of plane homogeneous Poisson structures.
  Inspired by   an idea  of I.~Vaisman \cite{Vaisman},   N.~Nakanishi \cite{Nakanishi} computed the cohomology of  plane quadratic Poisson structures classified by Z.-J.~Liu and P.~Xu \cite{LiuXu}.
   V.~I.~Arnold \cite{Arnold89} classified all plane Poisson structures with simple singularities in 1987 and Ph.~Monnier \cite{Monnier02a} computed the (germfied) Poisson cohomology of a large class of plane Poisson structures with isolated  singularities including Arnold's list. This work carries further  Ph.~Monnier's work  \cite{Monnier02a} by computing the Gerstenhaber algebra structure over the Poisson cohomology groups for these plane Poisson structures with isolated  singularities.

   \medskip

Let us give the setup of this paper.
Throughout,     $\mathbb{K}$ will  denote  the  field $\mathbb{R}$ or $\mathbb{C}$.
 We will consider the Poisson cohomology of a Poisson algebra $\mathcal{F}$, where
  $\mathcal{F}$, unless otherwise stated,   will always denote one of    the following algebras of  function germs:
\begin{itemize}\label{list of algebras}

\item $\mathbb{K}\{x, y\}$  the algebra of analytic functions near the origin,

\item $\mathbb{K}[[x, y]]$  the algebra of formal power series,

\item (only when $\mathbb{K}=\mathbb{R}$) $C_o^\infty(\mathbb{R}^2)$ the  algebra  of smooth function  germs at the origin.

\end{itemize}

Let  $\omega_{1},\omega_{2}$ be two positive integers, $W=\omega_{1}x\frac{\partial}{\partial x}+\omega_{2}y\frac{\partial}{\partial y}$ the Euler vector field. If a polynomial  $f$  satisfies  $Wf=\rmd f, $  with $\rmd\in\mathbb{Z}$, then $f$ is  \textit{weight-homogeneous}    of degree $\rmd$ with respect to $(\omega_{1},\omega_{2})$.
Denote
by
$\mathcal{P}_{\rmd-\omega_1-\omega_2}$ the space of weight-homogenous polynomials of degree $\rmd-\omega_1-\omega_2$ together with the zero polynomial,
and    write   $H_f=\frac{\partial f}{\partial y}\frac{\partial }{\partial x}-\frac{\partial f}{\partial x}\frac{\partial }{\partial y}$ for the Hamiltonian derivation associated to $f\in \mathcal{F}$.
  Recall that a polynomial $f $ is of finite codimension if the Milnor algebra $\mathcal{M}_f:= \mathcal{F}/I_f$ with $I_f=(\frac{\partial f}{\partial x}, \frac{\partial f}{\partial y})$ is finite-dimensional. Note that this implies that the origin is an isolated singularity.

\textit{From now on, we fix $(\omega_{1},\omega_{2})$ and all weight-homogenous stuff will be taken with respect to $(\omega_{1},\omega_{2})$.}

Let  $f$ be a weight-homogenous polynomial     of degree $\mathrm{d}$ (with respect to $(\omega_{1},\omega_{2})$).  Let  $h$ be a  weight-homogenous polynomial of degree $\mathrm{d}-\omega_{1}-\omega_{2}$ or the zero polynomial.  Suppose, moreover,  that $f$ has   finite codimension.

We shall consider the  Poisson structure   $$\Pi=f(1+h) \frac{\partial}{\partial x}\wedge\frac{\partial}{\partial y} $$ and if $h=0$, write  $\Pi_0:=\Pi=f   \frac{\partial}{\partial x}\wedge\frac{\partial}{\partial y}$. Our goal is to study  the Poisson cohomology   of this Poisson structure $\Pi$  for the algebras in the  above list.
For $i\geq 0$, let us denote the $i$-th Poisson cohomology group by $\PH^i_\Pi(\mathcal{F})$ or just  $\PH^i_\Pi $.
 Ph.~Monnier computed these Poisson cohomology groups.
\begin{Thm}[Ph.~Monnier \cite{Monnier02a}]
We have
 $$\left\{\begin{array}{rcl} \PH_{\Pi}^{0}&=&\k 1,\\
  \PH_{\Pi}^{1}&=&\k (1+h)H_f\oplus \mathcal{P}_{\rmd-\omega_1-\omega_2}(1+h)W,\\
   \PH_{\Pi}^{2}&=&\mathcal{M}_f\frac{\partial }{\partial x}\wedge \frac{\partial }{\partial y} \oplus \mathcal{P}_{\rmd-\omega_1-\omega_2}f \frac{\partial }{\partial x}\wedge \frac{\partial }{\partial y},\\
   \PH^i_\Pi&=&0, \ \mathrm{for}\  i\geq 3.\end{array}\right.$$
\end{Thm}
Therefore,  he actually showed that the dimensions of Poisson cohomology groups do not depend on the underlying algebras of function germs. This is a kind of GAGA type results \cite{Serre}.

\medskip

Poisson cohomology groups are not just abelian groups,  in fact they have a much nicer structure:  the so-called Gerstenhaber algebra structure \cite{Koszul,Xu99,LPV};  that is,   it   is a graded
commutative algebra
via the wedge product,    it has the Schouten-Nijenhuis  bracket
 so that it becomes  a graded Lie
algebra, and these structures are compatible.

A natural question is whether Ph.~Monnier's GAGA type result actually gives isomorphisms of Gerstenhaber algebras. The answer is Yes, which is the content of   the following result:
\begin{Thm}[Theorem~\ref{Thm: GAGA}]
 There exist isomorphisms of Gerstenhaber  algebras
$$ \PH_{\Pi}^{*}(\mathbb{C}\{x,y\})\cong \PH_{\Pi}^{*}(\mathbb{C}[[x,y]])\ \mathrm{and}\
 \PH_{\Pi}^{*}(C_{0}^{\infty}(\mathbb{R}^{2})) \cong \PH_{\Pi}^{*}(\mathbb{R}\{x,y\})\cong \PH_{\Pi}^{*}(\mathbb{R}[[x,y]]).$$
If $h=0$, one can also add the polynomial algebra  $\mathbb{K}[x, y]$   into the statement.
\end{Thm}

We prove the above result by  checking the difficult computation  of Ph.~Monnier \cite{Monnier02a} carefully. Moreover,
we  are able to  describe almost completely the Gerstenhaber algebra structure over the Poisson cohomology ring   $\PH_{\Pi}^{*}$.
\begin{Thm}
\begin{itemize}

\item[(a)](Proposition~\ref{Prop: Wedge product})  The wedge products $\wedge : \PH_{\Pi}^{0}\times \PH_{\Pi}^{i}\to \PH_{\Pi}^{i}, i=0, 1, 2$ are just  the scalar product of  $\k$-vector spaces;
the map $\wedge : \PH_{\Pi}^{1}\times \PH_{\Pi}^{1}\to \PH_{\Pi}^{2}$ can  be nonzero,   but its  only nonzero component   is given by
$$
 \mathbb{K}(1+h) H_{f}\times\mathcal{P}_{\rmd- \omega_{1}-\omega_{2}}(1+h)W \rightarrow \mathcal{P}_{\rmd- \omega_{1}-\omega_{2}}f\frac{\partial}{\partial x}\wedge\frac{\partial}{\partial y},\ \  ((1+h)H_f, g(1+h)W)\mapsto \rmd\cdot   gf \frac{\partial}{\partial x}\wedge\frac{\partial}{\partial y}$$
  for $g\in \mathcal{P}_{\rmd- \omega_{1}-\omega_{2}}$;
  all other wedge products   vanish.

  \item[(b)](Theorem~\ref{Thm: wedge product result}) We have an isomorphism of graded algebras $$\PH^{*}_{\Pi}\cong\k\langle u, v_1\cdots,v_{r}\rangle /(u^{2}, v_1^2, \cdots, v_{r}^{2},u v_{i}+v_{i}u, v_{i}v_{j},1\leq i,j\leq r)\times_{\k}\k[w_{1}]/(w_{1}^{2})\times_{\k}\cdots\times_{\k}\k[w_{c}]/(w_{c}^{2}),$$ with  $|u|=|v_1|=\cdots=|v_{r}|=1, |w_{1}|=\cdots=|w_{c}|=2$, where
      $r=\mathrm{dim}(\mathcal{M}_f), c=\mathrm{dim}(\mathcal{P}_{\rmd-\omega _1-\omega_2})$.

\item[(c)](Subsection~\ref{subsection: bracker for Pi_0}) If $h=0$,  then  all brackets vanish, except possibly the component  $$ [-,-]_{\sn}: \mathcal{P}_{\rmd-\omega _1-\omega_2} W \times  \mathcal{M}_f\frac{\partial }{\partial x}\wedge \frac{\partial }{\partial y}  \to \mathcal{M}_f\frac{\partial }{\partial x}\wedge \frac{\partial }{\partial y}.    $$
 of   $ [-,-]_{\sn}: \PH_{\Pi_0}^{1}\times \PH_{\Pi_0}^{2}\to \PH_{\Pi_0}^{2}$.

   \item[(d)](Remark~\ref{Rem: bracket for Pi}) If $h\neq 0$, then    all brackets vanish, except possibly the component of  $ [-,-]_{\sn}: \PH_{\Pi}^{1}\times \PH_{\Pi}^{2}\to \PH_{\Pi}^{2}$  given by    $$ [-,-]_{\sn}: \mathcal{P}_{\rmd-\omega _1-\omega_2}(1+h)W \times  \mathcal{M}_f\frac{\partial }{\partial x}\wedge \frac{\partial }{\partial y}  \to \mathcal{M}_f\frac{\partial }{\partial x}\wedge \frac{\partial }{\partial y} \oplus \mathcal{P}_{\rmd-\omega _1-\omega_2}f \frac{\partial }{\partial x}\wedge \frac{\partial }{\partial y}. $$ 
   \end{itemize}

\end{Thm}

When the isolated singularity is a simple singularity, we succeed in determining completely the Gerstenhaber algebra structure over the Poisson cohomology ring   $\PH_{\Pi}^{*}$ in terms of generators and relations; see  Table~\ref{Simple singularity computations} in Section~\ref{Section: Final result for simple singularities}.  Our calculation relies heavily on  the classification of plane Poisson structures with simple singularities due to V.~I.~Arnold  \cite{Arnold89}

\bigskip

This paper is organised as follows. The second section contains some basic definitions and facts about Poisson algebras and Poisson cohomology and we also remind  the classification of plane Poisson structures with simple singularities due to V.~I.~Arnold.  The result of Ph.~Monnier about Poisson cohomology groups is recalled in the third section. We consider the Gertstenhaber algebra structure  continuing the computation of  Ph.~Monnier  in the fourth  section. The last section contains explicit results for simple singularities.

\section{Preliminaries}
 For  basic definitions and facts about Poisson cohomology, the reader is referred to  Ref.~\onlinecite{LPV}.

 A Poisson algebra is an associative $\k$-algebra $R$ endowed  with a Lie bracket    $ \pi=\{\, ,\, \}$   satisfying a   compatibility condition.
For $p\geq 0$, a skew-symmetric  multilinear map $P\in Hom_{\k}(\wedge_\k^p R,R)$ is a $p$-polyvector field  if it is a derivation in each argument.
Denote by  $\mathfrak{X}^p(R)$  the space  of   $p$-polyvector fields.

Let us introduce the wedge product and the Schouten-Nijenhuis bracket over polyvector fields. 
\begin{Def}
  Let   $F\in \mathfrak{X}^p(R)$ and  $G\in \mathfrak{X}^q(R)$ with $p, q\geq 0$. Define the wedge product
~$F\wedge G\in \mathfrak{X}^{p+q}(R)$ as follows:
$$\begin{array}{rl}&F\wedge G(a_1\wedge \cdots \wedge a_{p+q})
 = \sum_{\sigma\in S_{p, q}}\mathrm{sgn}(\sigma) F(a_{\sigma(1)}\wedge\cdots \wedge a_{\sigma(p)})
G(a_{\sigma(p+1)}\wedge\cdots \wedge a_{\sigma(p+q)}),\end{array}$$ where $a_1, \cdots, a_{p+q}\in R$ and
$S_{p, q}$ is the set of $(p, q)$-shuffles, that is,
$$S_{p, q}=\{\sigma\in S_{p+q}\ |\ \sigma(1)<\cdots<\sigma(p), \sigma(p+1)<\cdots<\sigma(p+q)\}.$$

 Define the Schouten-Nijenhuis bracket $[F, G]_{\mathrm{\mathbf{SN}}}\in \mathfrak{X}^{p+q-1}(R)$ as follows:
 for $a_1, \cdots, a_{p+q-1}\in R$,
 $$\begin{array}{rl}
[F, G]_{\mathrm{\mathbf{SN}}}(a_1\wedge \cdots \wedge a_{p+q-1})
=&\sum\limits_{\sigma\in S_{q, p-1}}\mathrm{sgn}(\sigma) F(G(a_{\sigma(1)}\wedge\cdots \wedge a_{\sigma(q)})\wedge
 a_{\sigma(q+1)}\wedge\cdots \wedge a_{\sigma(p+q)})\\
&-(-1)^{(p-1)(q-1)}\sum\limits_{\sigma\in S_{p, q-1}}\mathrm{sgn}(\sigma)
  G(F(a_{\sigma(1)}\wedge\cdots \wedge a_{\sigma(p)})
 \wedge a_{\sigma(p+1)}\wedge\cdots \wedge a_{\sigma(p+q)}).\end{array}$$

\end{Def}


Since $\pi$ is a Lie bracket,   $\delta^p=[-,\pi]: \mathfrak{X}^{p}\to \mathfrak{X}^{p+1}, p\geq 0$ squares to zero, so  $(\mathfrak{X}^*(R), \delta^*)$ becomes a cochain complex, called
   the Poisson cochain complex     of $R$.
 For $p\geq 0$, the $p$-th Poisson cohomology group of the Poisson algebra $R$ is defined to be $$\PH^p(R):=\mathrm{Ker} (\delta^p)/\mathrm{Im} (\delta^{p-1}).$$
We write $\PH^*(R):=\bigoplus_{p\in \mathbb{N}} \PH^p(R).$
 It is easy to see that the wedge product and the Schouten-Nijenhuis bracket descend to cohomology groups.
The following result is folklore; see for example, Ref.~\onlinecite{LPV} Proposition 4.9.
\begin{Thm}\label{Thm: Poisson is Gerstenhaber}
Let $(R, \pi)$ be a Poisson algebra. Then $(\PH^* (R), \wedge, [, ]_{\mathrm{\mathbf{SN}}})$  is a  Gerstenhaber algebra, that is,
\begin{itemize}

\item[(a)] $(\PH^* (R), \wedge)$ is a graded commutative graded algebra,

\item[(b)]  $(\PH^* (R)[1],   [, ]_{\sn})$ is a graded Lie algebra,

\item[(c)]  the above two structures are compatible, that is,
$$[F\wedge G, H]_\sn=[F, H]_\sn\wedge G+(-1)^{(r-1)p}F\wedge [G, H]_\sn,$$
for $F\in \PH^p(R)$ and $H\in \PH^r(R)$.

\end{itemize}

\end{Thm}


Next we  recall the classification of plane Poisson structures with simple  singularities due to V.I.~Arnold \cite{Arnold89}. For unexplained notations and results about singularity theory, we refer the reader to Ref.~\onlinecite{AGV}.



  For plane Poisson structures with  simple singularities at the origin, V.I.~Arnold  proved the following classification result. 
 \begin{Thm}[V.I.Arnold \cite{Arnold89}]\label{Thm: Arnold3}

Let  $f$ be a simple function germ at  the origin of the plane.  The Poisson structure germ  $\Pi=f\frac{\partial}{\partial x}\wedge\frac{\partial}{\partial y}$ is equivalent to, up to a multiplicative constant,   $g\frac{\partial}{\partial x}\wedge\frac{\partial}{\partial y}$, where  $g$  is given in the second column of  Table~\ref{Table: Arnold's classification}:

\begin{table}[h]
 \centering
 \begin{tabular}{|c|c|c|c|c|c|c|}
\hline
type &  polynomial  & $\mathrm{d}$ &$\omega_1$ & $\omega_2$   \\
\hline
$A_{2p}\quad p\ge1$&$ x^{2}+y^{2p+1}$ &$4p+2$ & $2p+1$& $2$  \\
\hline
$A_{2p-1}^{\pm}\quad p\ge1$ & $(x^{2}\pm y^{2p})(1+\lambda y^{p-1})$ & $2p$ & $p$ & $1$\\
\hline
$D_{2p}^{\pm}\quad p\ge2$ & $(x^{2}y\pm y^{2p-1})(1+\lambda x+\mu y^{p-1})$ & $2p-1$&$p-1$ &$1$ \\
\hline
$D_{2p+1}\quad p\ge2$ & $(x^{2}y+y^{2p})(1+\lambda x)$&  $4p$&$2p-1$ &$2$  \\
\hline
$E_{6}$ & $x^{3}+y^{4}$  &  $12$&$4$ &$3$\\
\hline
$E_{7}$ & $(x^{3}+xy^{3})(1+\lambda y^{2})$ & $9$&$3$ &$2$ \\
\hline
$E_{8}$ & $x^{3}+y^{5}$ & $15$&$5$ &$3$ \\
\hline
\end{tabular}
\caption{Arnold's classification for plane Poisson structures with simple singularities}
\label{Table: Arnold's classification}
 \end{table}
Note that in this table, when  $\mathbb{K}=\mathbb{C}$, the symbol $\pm$ disappears.

\end{Thm}

\section{Poisson cohomology groups}

This section contains  a summary of the paper of Ph.~Monnier \cite{Monnier02a} and  we also explain how to compute cohomology classes concretely.

We want to investigate    the Poisson cohomology groups of $\Pi=f(1+h)\frac{\partial}{\partial x}\wedge\frac{\partial}{\partial y}$ and $\Pi_0=f\frac{\partial}{\partial x}\wedge\frac{\partial}{\partial y}$
where   $f$   is  a weight-homogeneous polynomial  of degree $\mathrm{d}>0$ which is a germ at the origin of finite codimension $\rm{c}$ and  $h$ is a weight-homogeneous polynomial  of degree $\mathrm{d}-\omega_{1}-\omega_{2}$ (only when  $\mathrm{d}-\omega_{1}-\omega_{2}>0$) or the zero polynomial. By Arnold's classification Theorem~\ref{Thm: Arnold3}, we see that the ``most interesting'' Poisson structures are of this type.

We choose a priori a monomial basis of $\mathcal{M}_f=\mathcal{F}/I_f$, say, $u_1, \cdots, u_c$ with $c=\mathrm{dim}(\mathcal{M}_f)$. For the existence of such a monomial basis, we refer the reader to Ref.~\onlinecite{AGV}. Note $\mathrm{max}(\mathrm{deg}(u_i), 1\leq i\leq  c)=2(\mathrm{d}-\omega_1-\omega_2)$. The choice of a monomial basis of $\mathcal{M}_f$ gives in fact a direct sum decomposition
$\mathcal{F}=\mathcal{M}_f\oplus I_f$. We choose a basis $e_1, \cdots, e_r$ of $\mathcal{P}_{\mathrm{d}-\omega_1-\omega_2}$ with $r=\mathrm{dim}(\mathcal{P}_{\mathrm{d}-\omega_1-\omega_2})$.

It is easy to see that  the space  $\PH^{0}_{\Pi_{0}}$ and $\PH^{0}_{\Pi}$ are always equal to $\mathbb{K}\cdot 1$.

The following  useful lemma is implicit in Ref.~\onlinecite{Monnier02a}, whose easy proof is left to the reader.
\begin{Lem}\label{Lem: intermidiate step}
Let $X$ be a weight-homogeneous vector field in the sense that $[W, X]_{\sn}=rX$ for some $r\in \mathbb{Z}$.
\begin{itemize}\item[(a)] Let $Z=\frac{1}{\mathrm{d}} (X(h)-\mathrm{div}(X)h)W+hX.$  Then
$\mathrm{div}(Z)=\frac{\deg(X)}{\mathrm{d}} (X(h)-\mathrm{div}(X))+2X(h)$ and
$X(fh)-\mathrm{div}(X)fh=Z(f).$

\item[(b)]  If  $\deg(X)\neq  \mathrm{d}-\omega_1-\omega_2$, denote $Y=X+ \frac{\mathrm{div}(X)}{\mathrm{d}-\omega_1-\omega_2- \mathrm{deg}(X)} W,$ then $X(f)=Y(f)- \mathrm{div}(Y)f.$

\item[(c)]
  If  $\deg(X)\neq 0$,    write $Y=hX-\frac{2}{\deg(X)}X(h)W,$ then  $\deg(Y)=\deg(X)+\mathrm{d}-\omega_1-\omega_2$ and  $$X(fh)-\mathrm{div}(X)fh= Y(f)- \mathrm{div}(Y)f.$$
\end{itemize}
\end{Lem}

With the same proof,  Theorems 4.5 and 4.9 of  Ref.~\onlinecite{Monnier02a}    hold  even for  $\mathbb{K}[x, y]$,  the polynomial ring in two variables.

\begin{Thm}[Ref.~\onlinecite{Monnier02a} Theorems 4.5 and 4.9 ]  \label{thm: H2 of Pi0} Let  $\mathcal{F}$ be either   $\mathbb{K}[x, y]$ or one of the algebras listed in the Introduction.
We have $$\PH^{1}_{\Pi_{0}}=\mathbb{K}H_f \oplus \mathcal{P}_{\rmd-\omega _1-\omega_2}W \ \mathrm{and}\
 \PH^{2}_{\Pi_{0}}=\mathcal{P}_{\rmd-\omega _1-\omega_2}f\frac{\partial}{\partial x}\wedge\frac{\partial}{\partial y}\oplus \mathcal{M}_f\frac{\partial}{\partial x}\wedge\frac{\partial}{\partial y}.$$


\end{Thm}

Let us include some details about how to compute the class of  $g\px\wedge \py$ in $\PH^2_{\Pi_0}$  for a weight-homogeneous polynomial $g\in \mathcal{F}$.

\begin{Prop}\label{Prop:Computing HP^2 of Pi0} Let $g\in \mathcal{F}$ be a weight-homogeneous polynomial.
Write $g=\sum_{i=1}^c \lambda_i u_i+\xi\in \mathcal{M}_f\oplus I_f=\mathcal{F}$ with $\lambda_i\in \k$   where $ \xi=X(f)\in I_f$ for a vector field $X$.
\begin{itemize}
\item[(a)]
If    $\mathrm{deg}(g)\neq 2\rmd-\omega_1-\omega_2$, then     $g\frac{\partial}{\partial x}\wedge\frac{\partial}{\partial y}=\sum_{i=1}^c \lambda_i u_i\frac{\partial}{\partial x}\wedge\frac{\partial}{\partial y}\in \PH^2_{\Pi_0}.$
\item[(b)] If $\mathrm{deg}(\xi)=2\rmd-\omega_1-\omega_2$, then   $g\frac{\partial}{\partial x}\wedge\frac{\partial}{\partial y}= \mathrm{div}(X)f\frac{\partial}{\partial x}\wedge\frac{\partial}{\partial y}\in \PH^2_{\Pi_0}.$\end{itemize}
\end{Prop}

\begin{Proof}
 Observe three facts:   if $\mathrm{deg}(g)> 2(\rmd-\omega_1-\omega_2)=\mathrm{max}(u_i, 1\leq i\leq c)$, then all $\lambda_i$ vanish; if $\mathrm{deg}(g)\leq 2(\rmd-\omega_1-\omega_2)$,
 those $u_i$ with nonzero $\lambda_i$ are of degree equal to that of $g$; the degree of $\xi$ is equal to $\mathrm{deg}(g)$.

  If  $\mathrm{deg}(\xi)\neq 2\rmd-\omega_1-\omega_2$, by Lemma 4.4 of Ref.~\onlinecite{Monnier02a}, $\xi\frac{\partial}{\partial x}\wedge\frac{\partial}{\partial y}$   is necessarily a $2$-coboundary.  By Lemma~\ref{Lem: intermidiate step}(b), let
$Y=X+\frac{\mathrm{div}(X)}{2\rmd-\omega_1-\omega_2- \mathrm{deg}(g)} W$, then
  $\xi=X(f)=Y(f)-\mathrm{div}(Y)f$.
Hence $$g\frac{\partial}{\partial x}\wedge\frac{\partial}{\partial y}=\sum_{i=1}^c \lambda_i u_i\frac{\partial}{\partial x}\wedge\frac{\partial}{\partial y}+(Y(f)-\mathrm{div}(Y)f)\frac{\partial}{\partial x}\wedge\frac{\partial}{\partial y}=\sum_{i=1}^c \lambda_i u_i\frac{\partial}{\partial x}\wedge\frac{\partial}{\partial y}+\delta^1_{\Pi_0}(f)=\sum_{i=1}^c \lambda_i u_i\frac{\partial}{\partial x}\wedge\frac{\partial}{\partial y}\in \PH^2_{\Pi_0}.$$

If $\mathrm{deg}(g)=\mathrm{deg}(\xi)=2\rmd-\omega_1-\omega_2>2(\rmd-\omega_1-\omega_2)$, then all $\lambda_i=0$. So
$g=\xi=\mathrm{div}(X)f+(X(f)-\mathrm{div}(X)f)$ and    $$g\frac{\partial}{\partial x}\wedge\frac{\partial}{\partial y}= \mathrm{div}(X)f\frac{\partial}{\partial x}\wedge\frac{\partial}{\partial y}\in \PH^2_{\Pi_0}.$$
\end{Proof}

Compared with the result for the Poisson structure $\Pi_0$,  computation of Ph.~Monnier  for  $\Pi$    cannot carry over $\mathbb{K}[x, y]$.
\begin{Thm}[Ref.~\onlinecite{Monnier02a}  Theorems 4.6 and 4.11]  \label{Thm: H2 of Pi}    We have $$\PH^{1}_{\Pi}=\mathbb{K}(1+h)H_f \oplus \mathcal{P}_{\rmd-\omega_1-\omega_2}(1+h)W\ \mathrm{and}\
 \PH^{2}_{\Pi}=\mathcal{P}_{\rmd-\omega_1-\omega_2}f\frac{\partial}{\partial x}\wedge\frac{\partial}{\partial y}\oplus \mathcal{M}_f\frac{\partial}{\partial x}\wedge\frac{\partial}{\partial y}.$$


\end{Thm}

We will explain how to compute the class of   $g\frac{\partial}{\partial x}\wedge\frac{\partial}{\partial y}    $ in $\PH^2_\Pi$ for $g$ a weight-homogeneous polynomial.
We need a small lemma which refines Lemma 4.10 of Ref.~\onlinecite{Monnier02a}.
\begin{Lem}\label{Lem:refinement} Let $Z$ be a weight-homogeneous vector field.
\begin{itemize}
\item[(a)] If $\deg(Z)>0$, then
    $ (Z(fh)-\mathrm{div}(Z)fh)\px\wedge\py=0\in \PH^2_{\Pi}$.
 \item[(b)] If $\deg(Z)=0$, then
  $ (Z(fh)-\mathrm{div}(Z)fh)\px\wedge\py =2Z(h)f\in \PH^2_{\Pi}.$
 \end{itemize}
\end{Lem}

\begin{Proof}
 First denote $g'=Z(fh)-\mathrm{div}(Z) fh$.   By  Theorem~\ref{thm: H2 of Pi0},  write
$$\frac{g'}{1+h}\px\wedge \py=\sum_{i=1}^c \mu_i u_i\px\wedge \py +Qf \px\wedge \py\in \PH^2_{\Pi_0}$$
with $\mu_i\in \k, Q\in \mathcal{P}_{\rmd-\omega_1-\omega_2}$.
In fact, it is shown  in the proof of Lemma 4.10 of  Ref.~\onlinecite{Monnier02a} that all $\mu_i$ vanish.

Let $g''$ be the component of $\frac{g'}{1+h}$
of degree $2\rmd-\omega_1-\omega_2$, then  following Proposition~\ref{Prop:Computing HP^2 of Pi0},  $Q=\mathrm{div}(Z')$  with  $g''=Z'(f)$, where  by Lemma~\ref{Lem: intermidiate step}(a),
$Z'=\frac{1}{\rmd} (Z(h)-\mathrm{div}(Z)h)W+hZ.$

 Express
$\frac{g'}{1+h}=g'-g'h+g'h^2+\cdots.$
If $\deg(Z)>0$, then   $\mathrm{deg}(g')>2\rmd -\omega_1-\omega_2$ and $\frac{g'}{1+h}$
has all its components of degree strictly bigger than $2\rmd -\omega_1-\omega_2$, and we obtain that $g''=0$ and $Q=0$. Hence, $ g'  \px\wedge \py=0\in \PH^2_{\Pi}$.
If $\deg(Z)=0$, then   $\mathrm{deg}(g')=2\rmd -\omega_1-\omega_2$ and $\frac{g'}{1+h}$
has all its components of degree strictly bigger than $2\rmd -\omega_1-\omega_2$ except $g'$, and we obtain that $g''=g'$ and $Q=\mathrm{div}(Z')=2Z(h)$. Hence, $ g'  \px\wedge \py=2Z(h)f\in \PH^2_{\Pi}.$

\end{Proof}

The following result  is very useful in   practical computation.
\begin{Prop} \label{Prop: final computation} Let $g\in \mathcal{F}$ be a weight-homogeneous polynomial.  Write $g=\sum_{i=1}^c \lambda_i u_i+\xi\in \mathcal{M}_f\oplus I_f=\mathcal{F}$ with   $\lambda_{i}\in \k$ and $\xi=X(f)\in I_f$ for a vector field $X$.

If  $\deg(g)$ is not in the form of $\rmd-k(\rmd-\omega_1-\omega_2)$ with $k\geq -1$, then  $$g\frac{\partial}{\partial x}\wedge\frac{\partial}{\partial y} =\sum_{i=1}^c \lambda_i u_i\frac{\partial}{\partial x}\wedge\frac{\partial}{\partial y}\in \PH^2_{\Pi}.$$

Suppose that  there exists an integer $k\geq -1$ such that $\deg(g)\rm=d-k(\rmd-\omega_1-\omega_2)$.
\begin{itemize}\item
  If $k=-1$, i.e.  $\deg(g)=2\rmd-\omega_1-\omega_2$, then
$$g\frac{\partial}{\partial x}\wedge\frac{\partial}{\partial y}=  \mathrm{div}(X)f\frac{\partial}{\partial x}\wedge\frac{\partial}{\partial y} \in \PH^2_{\Pi}. $$
 \item If $k>-1$, define   vector fields $X_i, 0\leq i\leq k$  as follows:
   put $X_0=X+\frac{\mathrm{div}(X)}{\rmd-\omega_1-\omega_2- \mathrm{deg}(X)} W$
   and  for $i=1, \cdots, k$,
 define recursively $X_{i}=hX_{i-1}-\frac{2}{\deg(X_{i-1})}X_{i-1}(h)W.$  Then $$g\frac{\partial}{\partial x}\wedge\frac{\partial}{\partial y}= \sum_{i=1}^c \lambda_i u_i\frac{\partial}{\partial x}\wedge\frac{\partial}{\partial y} -2X_k(h)f\frac{\partial}{\partial x}\wedge\frac{\partial}{\partial y} \in \PH^2_{\Pi}. $$
  \end{itemize}

\end{Prop}

\begin{Proof}


  Suppose that  $\deg(g)=\deg(\xi)=2\rmd-\omega_1-\omega_2=\rmd-k(\rmd-\omega_1-\omega_2)$ with $k=-1$. As  $$\begin{array}{rcl}\delta^1_{\Pi}(X)
  =(X(f)-\mathrm{div}(X)f)\px\wedge \py+(X(fh)-\mathrm{div}(X)fh)\px\wedge \py,\end{array}$$
  we have
 $$\begin{array}{rcl}\xi\px\wedge \py
 =\mathrm{div}(X)f\px\wedge \py-(X(fh)-\mathrm{div}(X)fh)\px\wedge \py\in \PH^2_{\Pi}.\end{array}$$
  Since $\deg(X)=\rmd-\omega_1-\omega_2>0$, by Lemma~\ref{Lem:refinement} (a), $(X(fh)-\mathrm{div}(X)fh)\px\wedge \py=0\in \PH^2_{\Pi}.$

  Suppose that $\deg(g)=\deg(\xi)\neq 2\rmd-\omega_1-\omega_2.$  Let
  $X_0=X+\frac{\mathrm{div}(X)}{2\rmd-\omega_1-\omega_2- \mathrm{deg}(g)} W.$
 Then by Lemma~\ref{Lem: intermidiate step}(b),  $\xi=X(f)=X_0(f)-\mathrm{div}(X_0)f$ and as in the previous paragraph,
    we need to consider the class of $(X_0(fh)-\mathrm{div}(X_0)fh)\px\wedge \py$ in $\PH^2_\Pi$.

 By Lemma~\ref{Lem:refinement} (a),  if  $\deg(X_0)>0$ or
  $\deg(g)=\deg(\xi)> 2\rmd-\omega_1-\omega_2=d-(-1)\cdot (\rm\rmd- \omega_1-\omega_2),$ we have
   $$(X_0(fh)-\mathrm{div}(X_0)fh)\px\wedge \py=0\in \PH^2_{\Pi} \  \mathrm{and}\
   \xi\px\wedge \py  =0 \in \PH^2_{\Pi}.$$

  If $\deg(X_0)=0$, that is $\deg(g)=\deg(\xi)= \rmd- 0\cdot (\rmd-\omega_1-\omega_2),$
  then by  Lemma~\ref{Lem:refinement} (b),
  $$(X_0(fh)-\mathrm{div}(X_0)fh)\px\wedge \py=2X_0(h)f\px\wedge \py\in \PH^2_{\Pi}.$$ Hence,
$\xi\px\wedge \py  =-2X_0(h)f \px\wedge \py  \in\PH^2_{\Pi}.$

If $\deg(X_0)<0$, that is,  $\deg(g)=\deg(\xi)< \rmd- 0\cdot (\rmd-\omega_1-\omega_2)$,
by Lemma~\ref{Lem: intermidiate step} (c), define $X_1=hX_0-\frac{2}{\deg(X_0)}X_0(h)W,$ then  $\deg(X_1)=\deg(X_0)+\rmd-\omega_1-\omega_2$ and  $X_0(fh)-\mathrm{div}(X_0)fh= X_1(f)- \mathrm{div}(X_1)f.$
 As above,  by Lemma~\ref{Lem:refinement} (a), if $\deg(X_1)>0$, that is, $ \rmd- 0\cdot (\rmd-\omega_1-\omega_2)>\deg(g)> \rmd- 1\cdot (\rmd-\omega_1-\omega_2)$, we obtain that $$(X_1(fh)-\mathrm{div}(X_1)fh)\px\wedge \py=0\in \PH^2_{\Pi}, \ \mathrm{and}\
  \xi\px\wedge \py  =0    \in\PH^2_{\Pi}.$$

 If $\deg(X_1)=0$,  by Lemma~\ref{Lem:refinement} (b),  $(X_1(fh)-\mathrm{div}(X_1)fh)\px\wedge \py=2X_1(h)\in \PH^2_{\Pi}$ and
 $\xi\px\wedge \py  =-2X_1(h)f\in \PH^2_{\Pi}.$
 If $\deg(X_1)<0$, i.e. $ \deg(g)< \rmd- 1\cdot (\rmd-\omega_1-\omega_2)$, by Lemma~\ref{Lem: intermidiate step} (c),  define $X_2=hX_1-\frac{2}{\deg(X_2)}X_2(h)W$ and we can continue as above.

 Let $k$ be the first integer  such that
$\mathrm{deg}(X_k)\geq 0$.
Then for $0\leq i\leq k-1$, $\deg(X_i)<0$, we see that
$$\xi\frac{\partial}{\partial x}\wedge\frac{\partial}{\partial y}=  -(X_k(fh)-\mathrm{div}(X_k) fh))\frac{\partial}{\partial x}\wedge\frac{\partial}{\partial y}  \in \PH^2_\Pi. $$
 We distinguish two cases. When  $\mathrm{deg}(X_k)> 0$, by Lemma~\ref{Lem:refinement} (a),
 $\xi\frac{\partial}{\partial x}\wedge\frac{\partial}{\partial y}=  -(X_k(fh)-\mathrm{div}(X_k) fh))\frac{\partial}{\partial x}\wedge\frac{\partial}{\partial y}=0  \in \PH^2_\Pi. $
 If $\mathrm{deg}(X_k)=0$, i.e. $\deg(g)=\rmd- k\cdot (\rmd-\omega_1-\omega_2)$,  by Lemma~\ref{Lem:refinement} (b),
 $ (X_k(fh)-\mathrm{div}(X_k) fh))\frac{\partial}{\partial x}\wedge\frac{\partial}{\partial y}= 2X_k(h)  \in \PH^2_\Pi, $
and
 $$\xi\frac{\partial}{\partial x}\wedge\frac{\partial}{\partial y}=  -(X_k(fh)-\mathrm{div}(X_k) fh))\frac{\partial}{\partial x}\wedge\frac{\partial}{\partial y}=-2X_k(h) f\px\wedge \py  \in \PH^2_\Pi. $$

\end{Proof}


\section{Wedge products  and Schouten-Nijenhuis   brackets  on Poisson cohomology}

\subsection{Wedge products}

 \begin{Prop} \label{Prop: Wedge product}

The wedge products $\wedge : \PH_{\Pi}^{0}\times \PH_{\Pi}^{i}\to \PH_{\Pi}^{i}, i=0, 1, 2$ are just  the scalar product of  $\k$-vector spaces;
the map $\wedge : \PH_{\Pi}^{1}\times \PH_{\Pi}^{1}\to \PH_{\Pi}^{2}$ can  be nonzero,   but its  only nonzero component   is given by
$$
 \mathbb{K}(1+h) H_{f}\times\mathcal{P}_{\rmd- \omega_{1}-\omega_{2}}(1+h)W \rightarrow \mathcal{P}_{\rmd- \omega_{1}-\omega_{2}}f\frac{\partial}{\partial x}\wedge\frac{\partial}{\partial y},\ \  ((1+h)H_f, g(1+h)W)\mapsto \rmd\cdot   gf \frac{\partial}{\partial x}\wedge\frac{\partial}{\partial y}$$
  for $g\in \mathcal{P}_{\rmd- \omega_{1}-\omega_{2}}$;
  all other wedge products   vanish.

\end{Prop}

\begin{Proof}
Since $\PH_{\Pi}^{0}=\k 1$,  the maps $\wedge : \PH_{\Pi}^{0}\times \PH_{\Pi}^{i}\to \PH_{\Pi}^{i}, i=0, 1, 2$ are just  the scalar product of  $\k$-vector spaces.

  It is easy to see that
$ (1+h)H_{f}\wedge(1+h)gW
= \rmd(1+h)^{2}gf\frac{\partial}{\partial x}\wedge\frac{\partial}{\partial y}
 .$
Since $\mathrm{deg}(hgf)=3\rmd- 2\omega_{1}-2\omega_{2}$, $\mathrm{deg}(h^{2}gf)=4\rmd- 3\omega_{1}-3\omega_{2}$ can not be written as $\rmd- k\cdot ( \rmd- \omega_{1}-\omega_{2})$ with $k\geq -1$ and are both greater than $2(\rmd-\omega_1-\omega_2)$,
by Proposition~\ref{Prop: final computation},
$\rmd(2h+h^2)gf\frac{\partial}{\partial x}\wedge\frac{\partial}{\partial y}=0\in \PH^2_\Pi.$
We obtain
$\rmd(1+h)^{2}gf\frac{\partial}{\partial x}\wedge\frac{\partial}{\partial y} =\rmd gf\frac{\partial}{\partial x}\wedge\frac{\partial}{\partial y}.$

The other statements are easy to show.
\end{Proof}

We can determine the algebra structure over  $\PH^{*}_{\Pi}$ in terms of generators and relations.  To this end,  we need a notation. Given two augmented $\k$-algebra $A$ and $B$, i.e. $\k$-algebras equipped with a $\k$-algebra homomorphism to $\k$,   the algebra
$A\times_\k B$ is the fiber product in the category of augmented $\k$-algebras.
$$\xymatrix{ A\times_\k B\ar[d]\ar[r]& B\ar[d] \\ A\ar[r]& \k}$$

\begin{Thm} \label{Thm: wedge product result} We have an isomorphism of graded algebras $$\PH^{*}_{\Pi}\cong\k\langle u, v_1\cdots,v_{r}\rangle /(u^{2}, v_1^2, \cdots, v_{r}^{2},u v_{i}+v_{i}u, v_{i}v_{j},1\leq i,j\leq r)\times_{\k}\k[w_{1}]/(w_{1}^{2})\times_{\k}\cdots\times_{\k}\k[w_{c}]/(w_{c}^{2}),$$ with  $|u|=|v_1|=\cdots=|v_{r}|=1, |w_{1}|=\cdots=|w_{c}|=2$. Here $|u|$ denotes the cohomological degree of $u$.

\end{Thm}

\begin{Proof}
Denote $u= (1+h)H_{f}, v_{i}=(1+h)e_{i}W, 1\le i\le r,  w_{j}=  u_{j}\frac{\partial}{\partial x}\wedge\frac{\partial}{\partial y},1\le j\le c$.  By the preceding proposition, we obtain the result.
\end{Proof}



\subsection{Schouten-Nijenhuis brackets  on $\PH^{*}_{\Pi_{0}}$} \label{subsection: bracker for Pi_0}
By dimension reason and as $\PH_{\Pi_{0}}^{0}=\mathbb{K}$, we only need to compute
$[\ ,\ ]_{\sn}:  \PH_{\Pi_{0}}^{1}\times  \PH_{\Pi_{0}}^{1}\to  \PH_{\Pi_{0}}^{1}$ and
$[\ ,\ ]_{\sn}:  \PH_{\Pi_{0}}^{1}\times  \PH_{\Pi_{0}}^{2}\to  \PH_{\Pi_{0}}^{2}.$
We need two formulae in the sequel:
 for all $f,g\in\mathcal{F}$,
 \begin{equation}\label{eq: []}
 [f\frac{\partial}{\partial x}, g\frac{\partial}{\partial x}\wedge\frac{\partial}{\partial y}]_{\sn}=(f\frac{\partial g}{\partial x}-g\frac{\partial f}{\partial x})\frac{\partial}{\partial x}\wedge\frac{\partial}{\partial y}, \ \ \  [f\frac{\partial}{\partial y}, g\frac{\partial}{\partial x}\wedge\frac{\partial}{\partial y}]_{\sn}=(f\frac{\partial g}{\partial y}-g\frac{\partial f}{\partial y})\frac{\partial}{\partial x}\wedge\frac{\partial}{\partial y}.
\end{equation}

\begin{Lem} The map $[\ ,\ ]_{\sn}:  \PH_{\Pi_{0}}^{1}\times  \PH_{\Pi_{0}}^{1}\to  \PH_{\Pi_{0}}^{1}$ vanishes.

\end{Lem}

\begin{Proof} We will deal with the most difficult case
$[H_{f},e_{i}W]_{\sn},~1\leq i\leq r
 $, the other cases being easy.

 We have
$ [H_{f},e_{i}W]_{\sn}=-\rmd fH_{e_{i}}.$
Since $\deg(\rmd fH_{e_{i}})=2\rmd- 2\omega_{1}-2\omega_{2}> \rmd- \omega_{1}-\omega_{2},$ by  \onlinecite{Monnier02a} Lemma 4.4, $\rmd fH_{e_{i}}=0\in \PH^1_{\Pi_0}$.
\end{Proof}

\begin{Prop} \label{Prop: bracket for Pi_0} All brackets in $[\ ,\ ]_{\sn}:  \PH_{\Pi_{0}}^{1}\times  \PH_{\Pi_{0}}^{2}\to  \PH_{\Pi_{0}}^{2}$ vanish except   possibly the component   $$ [-,-]_{\sn}: \mathcal{P}_{\rmd-\omega_1-\omega_2} W \times  \mathcal{M}_f\frac{\partial }{\partial x}\wedge \frac{\partial }{\partial y}  \to \mathcal{M}_f\frac{\partial }{\partial x}\wedge \frac{\partial }{\partial y}.   $$ 
\end{Prop}

\begin{Proof}
We need to show the following statements:
\begin{itemize}
\item[(a)]  For $1\leq i\leq c$,  $[u, w_i]_\sn=[H_{f},u_{i}\frac{\partial}{\partial x}\wedge\frac{\partial}{\partial y}]_{\sn}=0\in \PH^2_{\Pi_0}$;

\item[(b)] For $1\leq i\leq r$, $[H_{f}, e_{i}f\frac{\partial}{\partial x}\wedge\frac{\partial}{\partial y}]_{\sn}=0 \in \PH^2_{\Pi_0};$

    \item[(c)]  For $ 1\leq i,j\leq r$, $
[e_{i}W,e_{j}f\frac{\partial}{\partial x}\wedge\frac{\partial}{\partial y}]_{\sn}= 0 \in \PH^2_{\Pi_0}.$

    \end{itemize}
In fact, as $\mathrm{div}(H_{u_i})=0$,
$$[u, w_i]_\sn=[H_{f},u_{i}\frac{\partial}{\partial x}\wedge\frac{\partial}{\partial y}]_{\sn}= -H_{u_i}(f)\frac{\partial}{\partial x}\wedge \frac{\partial}{\partial y}= -(H_{u_i}(f)-\mathrm{div}(H_{u_i})f)\frac{\partial}{\partial x}\wedge \frac{\partial}{\partial y}=-\delta^1(H_{u_i}),$$
 so $[H_{f},u_{i}\frac{\partial}{\partial x}\wedge\frac{\partial}{\partial y}]_{\sn}=0\in \PH^2_{\Pi_0}$.
 The other two statements follow from $w_i= e_{i}f\frac{\partial}{\partial x}\wedge\frac{\partial}{\partial y}=\frac{1}{d}u\wedge v_i$ and Theorem~\ref{Thm: Poisson is Gerstenhaber} (c).

  Using (\ref{eq: []}), it is not difficult to see that $[e_{i}W, u_{j}\frac{\partial}{\partial x}\wedge\frac{\partial}{\partial y}]_{\sn}=(\mathrm{deg}(u_{j})-d)e_{i}u_{j}\frac{\partial}{\partial x}\wedge \frac{\partial}{\partial y}.$
We can write
$(\mathrm{deg}(u_{j})-d)e_{i}u_{j} =\sum_k \lambda_k u_k+\xi\in M_f\oplus I_f=\mathcal{F}$
with $\lambda_k\in \mathbb{K}$ and $\xi\in I_f$.
  If  $\rmd- \omega_{1}-\omega_{2}+\deg(u_{j})=deg(\xi)=2\rmd-\omega_1-\omega_2$ and  thus $\deg(u_{j})=d$, but then $(\deg(u_{j})-d)e_{i}u_{j} =0$ and so $\xi=0\in \PH^2_{\Pi_0}$.
If $\deg(u_j)\neq d$, by  Proposition~\ref{Prop:Computing HP^2 of Pi0},
$$[e_{i}W, u_{j}\frac{\partial}{\partial x}\wedge\frac{\partial}{\partial y}]_{\sn}=\sum_k \lambda_k u_k\in \PH^2_{\Pi_0}.$$ Of course only those $u_k$ with $\deg(u_k)=\deg(u_j)+\rmd- \omega_{1}-\omega_{2}$ could occur.


\end{Proof}

\subsection{Schouten-Nijenhuis brackets  on $\PH^{*}_\Pi$}\label{subsection bracket for Pi}

Remember that in this case, we always assume $\rmd-\omega_1-\omega_2>0$.

By dimension reason and as $\PH_{\Pi}^{0}=\mathbb{K}$, we only need to consider
$[\ ,\ ]_{\sn}:  \PH_{\Pi}^{1}\times  \PH_{\Pi}^{1}\to  \PH_{\Pi}^{1}$ and
$[\ ,\ ]_{\sn}:  \PH_{\Pi}^{1}\times  \PH_{\Pi}^{2}\to  \PH_{\Pi}^{2}.$

\begin{Lem} The map $[\ ,\ ]_{\sn}:  \PH_{\Pi}^{1}\times  \PH_{\Pi}^{1}\to  \PH_{\Pi}^{1}$ vanishes.

\end{Lem}

\begin{Proof} We will deal with the most difficult case
$[u, v_i]_\sn,~1\leq i\leq r
 ,$ the other cases being easy.

  A direct computation shows that $$[u, v_i]_\sn=
 [(1+h)H_{f},(1+h)e_{i}W]_{\sn}= -\rmd (1+h)fH_{e_{i}(1+h)}.$$
Since the degrees of $ \rmd fH_{e_{i}}$  and  $ \rmd fH_{e_{i}h}$     are all strictly larger than $ \rmd- \omega_{1}-\omega_{2}$, by Ref.~\onlinecite{Monnier02a} Lemma 4.4,  $\rmd fH_{e_{i}(1+h)}=0\in  \PH^1_{\Pi_0}$. By (the proof of) Theorem 4.6 of Ref.~\onlinecite{Monnier02a}, $$ \rmd(1+h)fH_{e_{i}(1+h)}\in B^1_\Pi  \Leftrightarrow
  \frac{\rmd (1+h)f}{1+h}H_{e_{i}(1+h)}=\rmd fH_{e_{i}(1+h)}\in B^1_{\Pi_0}.$$
    Hence $[(1+h)H_{f},(1+h)e_{i}W]_{\sn}=0\in \PH^1_{\Pi}.$
 \end{Proof}

\begin{Lem}We have
\begin{itemize}
\item[(a)] For $1\leq i\leq r$, $[u, e_{i}f\frac{\partial}{\partial x}\wedge\frac{\partial}{\partial y}]_{\sn}=[(1+h)H_{f}, e_{i}f\frac{\partial}{\partial x}\wedge\frac{\partial}{\partial y}]_{\sn}=0\in \PH^2_\Pi$.

\item[(b)] For $1\leq i,j\leq r$, $[v_i, e_{j}f\frac{\partial}{\partial x}\wedge\frac{\partial}{\partial y}]_{\sn}=[(1+h)e_{i}W, e_{j}f\frac{\partial}{\partial x}\wedge\frac{\partial}{\partial y}]_{\sn}=0\in \PH^2_\Pi.$\end{itemize}
 \end{Lem}
\begin{Proof}
 The statements follow from Theorem~\ref{Thm: Poisson is Gerstenhaber} (c) and the previous lemma, as $ e_{i}f\frac{\partial}{\partial x}\wedge\frac{\partial}{\partial y}=\frac{1}{d} u\wedge v_i$.
\end{Proof}

\begin{Lem}
For all $1\leq i\leq c$, $[u, w_i]_\sn=[(1+h)H_{f},u_{i}\frac{\partial}{\partial x}\wedge\frac{\partial}{\partial y}]_{\sn}=0\in \PH^2_\Pi$.
\end{Lem}

\begin{Proof}  Using (\ref{eq: []}), it is not difficult to see that
$$\begin{array}{rcl}
 [(1+h)H_{f},u_{i}\frac{\partial}{\partial x}\wedge\frac{\partial}{\partial y}]_{\sn} &
=&(-(1+h)H_{u_i}(f)+u_{i}H_{h}(f))\frac{\partial}{\partial x}\wedge \frac{\partial}{\partial y}\\
&=&-\delta^1(H_{u_i})+(-f H_{h}(u_i)+u_{i}H_{h}(f))\frac{\partial}{\partial x}\wedge \frac{\partial}{\partial y}\\
& =&X (f)\frac{\partial}{\partial x}\wedge \frac{\partial}{\partial y} \, \mathrm{mod}\, B^2_{\Pi},
\end{array}$$ where $X =-\frac{1}{d} H_{h}(u_i)W +u_{i}H_{h}$.
 Notice that $\mathrm{div}(X )= \frac{\omega_{1}+\omega_{2}-\mathrm{deg}(u_{i})}{d} H_h(u_i).$

As in Proposition~\ref{Prop: final computation}, denote $g=X(f)\in I_f$. Then
  $\mathrm{deg}(g )=2 \rmd- 2\omega_{1}-2\omega_{2} +\mathrm{deg}(u_{i})$.
  It is equal to $\rmd- k(\rmd-  \omega_{1}- \omega_{2})$ with $k\geq -1$ if and only if
  $\deg(u_i)=\rmd- (k+2)(\rmd-  \omega_{1}- \omega_{2})$.

  If $\mathrm{deg}(g )$
  is not of the form $\rmd- k(\rmd-  \omega_{1}- \omega_{2})$ with $k\geq -1$.
  By Proposition~\ref{Prop: final computation},   $g\frac{\partial}{\partial x}\wedge \frac{\partial}{\partial y}= 0\, \in \PH^2_\Pi .$

  Now suppose that $\mathrm{deg}(g )= \rmd- k(\rmd-  \omega_{1}- \omega_{2})$ with $k\geq -1$.

Let $k=-1$, i.e.
  $\mathrm{deg}u_{i}=\omega_{1}+\omega_{2}$ and $\mathrm{deg}(g )=2 \rmd- \omega_{1}- \omega_{2}$.
  Then $\mathrm{div(X)}=0$ and
   $$g\frac{\partial}{\partial x}\wedge \frac{\partial}{\partial y}=X (f)\frac{\partial}{\partial x}\wedge \frac{\partial}{\partial y}= \mathrm{div}(X)f\frac{\partial}{\partial x}\wedge \frac{\partial}{\partial y}= 0\, \in \PH^2_\Pi .$$


Let $k=0$. Then $\mathrm{deg}(u_i)= 2\omega_1+2\omega_2-d$. Let $X_0=X+\frac{\mathrm{div}(X)}{\rmd-\omega_1-\omega_2-\deg(X)}W= u_iH_h.$ Then
$X(f)=X_0(f)-\mathrm{div}(X_0)f$.   Remark that $X_0(h)=0$ and
  $X_0(fh)-\mathrm{div}(X_0)fh=hX(f). $ Hence,
 $$g\frac{\partial}{\partial x}\wedge \frac{\partial}{\partial y}=  -2X_0(h)f\frac{\partial}{\partial x}\wedge \frac{\partial}{\partial y}= 0\, \in \PH^2_\Pi .$$

Suppose that $k\geq 1$. Let $X_1=hX_{0}-\frac{2}{\deg(X_{0})}X_{0}(h)W=hX_0.$

 If $k=1$, then $\deg(u_i)=-2d+3\omega_1+3\omega_2$  and
 $$g\frac{\partial}{\partial x}\wedge \frac{\partial}{\partial y}=  -2X_1(h)f\frac{\partial}{\partial x}\wedge \frac{\partial}{\partial y}=  -2hX_0(h)f\frac{\partial}{\partial x}\wedge \frac{\partial}{\partial y}=0\, \in \PH^2_\Pi .$$

  By induction, $X_i=h^iX_0$ for $1\leq i\leq k$ and $$g\frac{\partial}{\partial x}\wedge \frac{\partial}{\partial y}=  -2X_k(h)f\frac{\partial}{\partial x}\wedge \frac{\partial}{\partial y}=  -2h^kX_0(h)f\frac{\partial}{\partial x}\wedge \frac{\partial}{\partial y}=0\, \in \PH^2_\Pi.$$

\end{Proof}

  \begin{Lem}\label{Lem: [v_i, w_j]}   For $1\leq i\leq r$ and $ 1\leq j\leq c,$ we have
$$
[v_i, w_j]=[(1+h)e_{i}W,u_{j}\frac{\partial}{\partial x}\wedge\frac{\partial}{\partial y}]_{\sn}=((\mathrm{deg}u_{j}-d)e_{i}u_{j}+(\mathrm{deg}u_{j}-2d+\omega_{1}+\omega_{2})he_{i}u_{j})\frac{\partial}{\partial x}\wedge \frac{\partial}{\partial y}.
 $$
By Proposition~\ref{Prop: final computation}, there exists an algorithm to compute its class in $\PH^2_{\Pi}$.
\end{Lem}

\begin{Proof}
Using (\ref{eq: []}), we have
$$\begin{array}{lll}
[(1+h)e_{i}W,u_{j}\frac{\partial}{\partial x}\wedge\frac{\partial}{\partial y}]_{\sn}
=((\mathrm{deg}u_{j}-d)e_{i}u_{j}+(\mathrm{deg}u_{j}-2d+\omega_{1}+\omega_{2})he_{i}u_{j})\frac{\partial}{\partial x}\wedge \frac{\partial}{\partial y}.
\end{array}$$
We will use   Proposition~\ref{Prop: final computation} to compute its cohomology class  in $\PH^2_{\Pi}.$

Let $g_1=(\mathrm{deg}u_{j}-d)e_{i}u_{j}$ and $g_2=(\mathrm{deg}u_{j}-2d+\omega_{1}+\omega_{2})he_{i}u_{j}$.
Notice that  $\mathrm{deg}(g_1)=\mathrm{d}-\omega_{1}-\omega_{2}+\mathrm{deg}(u_{j})$,   $\mathrm{deg}(g_2)=2(\mathrm{d}-\omega_{1}-\omega_{2})+\mathrm{deg}(u_{j})$,   and $\mathrm{deg}(u_{j})=\rmd- k(\rmd- \omega_{1}-\omega_{2})$ if and only if  $\mathrm{deg}(g_1)= \rmd- (k-1)(\rmd- \omega_{1}-\omega_{2})$ if and only if
 $\mathrm{deg}(g_2)= \rmd- (k-2)(\rmd- \omega_{1}-\omega_{2})$.

 Write $g_1=\sum_{s=1}^c \lambda_s u_s +Y(f)$  and  $g_2=\sum_{s=1}^c \mu_s u_s +Z(f)$ for $\lambda_s, \mu_s\in \k$ and two vector fields $Y, Z$.

 If   there  does not  exist  $k\geq 0$ such that
     $\mathrm{deg}(u_{j})=\rmd- k(\rmd- \omega_{1}-\omega_{2})$, then
     $$[v_1, w_j]_{\sn}=\sum_{s=1}^c (\lambda_s+\mu_s)w_s\in \PH^2_\Pi.$$

 Suppose  there   exists $k\geq 0$ such that
     $\mathrm{deg}(u_{j})=\rmd- k(\rmd- \omega_{1}-\omega_{2})$. Then $\mathrm{deg}(g_1)= \rmd- (k-1)(\rmd- \omega_{1}-\omega_{2})$ and
 $\mathrm{deg}(g_2)= \rmd- (k-2)(\rmd- \omega_{1}-\omega_{2})$. We distinguish three cases: $k=0, k=1, k\geq 2$.

 Let $k=0$. Then $k-1=-1, k-2=-2<-1$.  We see
  $$[v_1, w_j]_{\sn}=\sum_{s=1}^c  \mu_s w_s +\mathrm{div}(Y)f\frac{\partial}{\partial x}\wedge\frac{\partial}{\partial y}    \in \PH^2_\Pi.$$

  Let $k=1$. Then $k-1=0, k-2=-1$.  Define as in Proposition~\ref{Prop: final computation},
  $Y_0=Y+\frac{\mathrm{div}(Y)}{\rmd-  \mathrm{deg}(u_j)} W.$ We have
  $$[v_1, w_j]_{\sn}=\sum_{s=1}^c \lambda_s w_s +(-2Y_0(h) +\mathrm{div}(Z))f\frac{\partial}{\partial x}\wedge\frac{\partial}{\partial y} \in \PH^2_{\Pi}.$$

  Let $k\geq 2$. Then $k-1\geq 1, k-2\geq 0$.
     Define  $$Y_{t}=hY_{t-1}-\frac{2}{\deg(Y_{t-1})}Y_{t-1}(h)W, t=1, \cdots, k-1,$$
        $$Z_0=Z+\frac{\mathrm{div}(Z)}{\rmd-\omega_1-\omega_2- \mathrm{deg}(Z)} W=Z+\frac{\mathrm{div}(Z)}{\omega_1+\omega_2-  \mathrm{deg}(u_j)} W$$
   and
 $$Z_{t}=hZ_{t-1}-\frac{2}{\deg(Z_{t-1})}Z_{t-1}(h)W, t=1, \cdots, k-2.$$
  We obtain
$$[v_1, w_j]_{\sn}=\sum_{s=1}^c (\lambda_s+\mu_s)w_s -2(Y_{k-1}(h)+Z_{k-2}(h)) f\frac{\partial}{\partial x}\wedge\frac{\partial}{\partial y} \in \PH^2_{\Pi}.$$

\end{Proof}

   \begin{Rem} \label{Rem: bracket for Pi}
  All brackets vanish except
   $$ [-,-]_{\sn}: \mathcal{P}_{\rmd-\omega_1-\omega_2}(1+h)W \times  \mathcal{M}_f\frac{\partial }{\partial x}\wedge \frac{\partial }{\partial y}  \to \mathcal{M}_f\frac{\partial }{\partial x}\wedge \frac{\partial }{\partial y} \oplus \mathcal{P}_{\rmd-\omega_1-\omega_2}f \frac{\partial }{\partial x}\wedge \frac{\partial }{\partial y}.$$ which is given by (the proof of) Lemma~\ref{Lem: [v_i, w_j]}.
\end{Rem}

Since our computation  in this   section does not depend on the function algebra $\mathcal{F}$, we have obtained the following GAGA type result:
\begin{Thm}\label{Thm: GAGA}
 There exist isomorphisms of Gerstenhaber  algebras
$$ \PH_{\Pi}^{*}(\mathbb{C}\{x,y\})\cong \PH_{\Pi}^{*}(\mathbb{C}[[x,y]])\ \mathrm{and}\
 \PH_{\Pi}^{*}(C_{0}^{\infty}(\mathbb{R}^{2})) \cong \PH_{\Pi}^{*}(\mathbb{R}\{x,y\})\cong \PH_{\Pi}^{*}(\mathbb{R}[[x,y]]).$$
If $h=0$, one can also add the polynomial algebra in two variables into the statement.
\end{Thm}

\section{Simple singularities}\label{Section: simple singularities}

For simple singularities, following the classification of V.I.~Arnold \cite{Arnold89}, we can determine the Gerstenhaber algebra structure over the Poisson cohomology of plane Poisson structures with simple singularities.

\subsection{Type $A_{2p}$:$x^{2}+ y^{2p+1}, p\geq 1$}

In this type $A_{2p}$, $\Pi_0=f\px\wedge \py, f=x^2+y^{2p+1},
   \omega_{1}=2p+1, \omega_{2}=2, d=4p, \rmd-\omega_1-\omega_2=2p.$
 We have   $\mathcal{P}_{ \rmd-\omega_1-\omega_2}=\mathcal{P}_{2p}=\k y^{2p}$,  $r=1, e_1=y^{2p}$.
 We also  have  $\mathcal{M}_f\cong \k[y]/(y^{2p})$, so $c=2p$ and choose  $u_i=y^{i-1}, i=1, \cdots, 2p$.
Let $$u= H_f=  (2p+1)y^{2p}\frac{\partial}{\partial x}-2x\frac{\partial}{\partial y},
    v_1=e_1 W=y^{2p} ((2p+1)x\frac{\partial}{\partial x}+2y\frac{\partial}{\partial y})$$ which spans   $\mathcal{P}_{\rmd-\omega_1-\omega_2} W$,  and  $w_i=u_i \frac{\partial}{\partial x}\wedge\frac{\partial}{\partial y}=y^{i-1}\frac{\partial}{\partial x}\wedge\frac{\partial}{\partial y}$ with $1\leq i\leq 2p $, which form a basis of $\mathcal{M}_f \frac{\partial }{\partial x}\wedge \frac{\partial }{\partial y}$.

By Theorem~\ref{Thm: wedge product result}   and Proposition~\ref{Prop: bracket for Pi_0}, we have the following result:
\begin{Prop} Let  $\Pi_0= (x^2+y^{2p+1})\px\wedge \py $ be Type $A_{2p}, p\geq 1$. We have an isomorphism of Gerstenhaber algebras  $$\PH^{*}_{\Pi_0}\cong\wedge_\k(u, v_{1})\times_{\k}\k[w_{1}]/(w_{1}^{2})\times_{\k}\cdots\times_{\k}\k[w_{2p}]/(w_{2p}^{2}),$$
with  $|u|=1, |w_{1}|= \cdots=|w_{2p }|=2$,
where the brackets are all zero.
\end{Prop}

\subsection{Type $A_{2p-1}^{\pm}$: $(x^{2}\pm y^{2p})(1+\lambda y^{p-1}), p\ge 1$}

In this type $A_{2p-1}^{\pm}$, $\Pi=f(1+h)\px\wedge \py $ with $f=x^2\pm y^{2p}, h=\lambda y^{p-1},
   \omega_{1}=p, \omega_{2}=1, d=2p, \rmd-\omega_1-\omega_2=p-1.$
 We have   $\mathcal{P}_{ \rmd-\omega_1-\omega_2}=\mathcal{P}_{p-1}=\k y^{p-1}$,  $r=1, e_1=y^{p-1}$.
 We also  have  $\mathcal{M}_f\cong \k[y]/(y^{2p-1})$, so $c=2p-1$ and choose  $u_i=y^{i-1}, i=1, \cdots, 2p-1$.
As in Theorem~\ref{Thm: wedge product result}, let $$u=(1+h)H_f=(1+\lambda y^{p-1})(\pm 2py^{2p-1}\frac{\partial}{\partial x}-2x\frac{\partial}{\partial y}),
     v_1=e_1(1+h)W=y^{p-1}(1+\lambda y^{p-1})(px\frac{\partial}{\partial x}+y\frac{\partial}{\partial y})$$ and
    $w_i=u_i \frac{\partial}{\partial x}\wedge\frac{\partial}{\partial y}=y^{i-1}\frac{\partial}{\partial x}\wedge\frac{\partial}{\partial y}, 1\leq i\leq 2p-1.$
We have an isomorphism of graded algebras  $$\PH^{*}_{\Pi}\cong\wedge_\k(u, v_{1})\times_{\k}\k[w_{1}]/(w_{1}^{2})\times_{\k}\cdots\times_{\k}\k[w_{2p-1}]/(w_{2p-1}^{2}).$$
Notice that by Proposition~\ref{Prop: Wedge product}, $u\wedge v_1=2p e_1f\px\wedge \py\in \PH^2_{\Pi}.$

\medskip

Let us first consider the case $p=1$. In this case, $\rmd-\omega_1-\omega_2=0$ and we are   in the situation of $\Pi_0=(x^2\pm y^{2})\px\wedge \py$. By Subsection~\ref{subsection: bracker for Pi_0}, the only nonzero bracket is
$[v_1, w_1]_\sn=[W, \frac{\partial}{\partial x}\wedge\frac{\partial}{\partial y}]=-2 w_1.$

Suppose now that $p>1$. By Subsection~\ref{subsection bracket for Pi},   we need to consider
 $[v_1, w_j]_{\sn}, j=1, \cdots, 2p-1$.
By Lemma~\ref{Lem: [v_i, w_j]},  for $1\leq j\leq 2p-1$,
$
[v_1, w_j]_{\sn}= ((j-2p-1)y^{p+j-2}+(j-3p)\lambda y^{2p+j-3})\frac{\partial}{\partial x}\wedge \frac{\partial}{\partial y}.
 $  As in Lemma~\ref{Lem: [v_i, w_j]}, let $g_1=(j-2p-1)y^{p+j-2}, g_2=(j-3p)\lambda y^{2p+j-3}$.
Then $\deg(g_1)=p+j-2$ and $\deg(g_2)=2p+j-3$.
We consider $k\in \mathbb{Z}$ such that   $\mathrm{deg}(g_1)= \rmd- k(\rmd- \omega_{1}-\omega_{2})=2p-k(p-1)$, then
   $\mathrm{deg}(g_2)=\rmd- (k-1)(\rmd- \omega_{1}-\omega_{2})=2p-(k-1)(p-1)$  and
   $\mathrm{deg}(u_{j})=\rmd- (k+1)(\rmd- \omega_{1}-\omega_{2})=2p-(k+1)(p-1)$ .
The restrictions are 
 $0\leq \deg(u_j)\leq 2(\rmd-\omega_1-\omega_2)=2p-2$ and the degrees of $g_1, g_2$ are all nonnegative.
It is easy to see that  all the possibilities are  listed in Table~\ref{tab1}:
 \begin{table}[h]
 \centering
 \begin{tabular}{|c|c|c|c|c|c|}
\hline
   k & $\deg(g_1)$ &$\deg(g_2)$ & $\deg(u_j)$ & $j$ &  $p $\\
\hline

\hline
 $1$ & $4 $ &$5 $ & $3 $&$ 4$  &$2 $\\
\hline
 $2$ & $2$ &$3$ & $1$&$2$  &$2 $\\
\hline
$3$ & $1$ & $2$ & $0$ &$1$  & $2$\\
\hline
$0$ &  $6 $ &$8 $ & $4 $ & $ 5$  &$ 3$\\
\hline
$1$ & $4 $ &$6 $ & $2 $&$ 3$  &$ 3$\\
\hline
$2$ & $2$ &$4$ & $0$&$1$  &$ 3$\\

\hline

$0$ & $2p$& $ 3p-1$ & $p+1$ & $p+2$ & $> 3$ \\
\hline
$1$ &$p+1$ & $2p$ &$2$ & $3$ & $>3$\\
\hline
\end{tabular}
\caption{All possibilities in the computation}
\label{tab1}
 \end{table}


\textbf{Case $k=1, p=2$.} In this case, $j=3, u_3=y^2$,   $g_1=-2y^{3}, g_2=-3\lambda y^{4}$.
We have $ g_2=Z(f)\in I_f$.
 with $ Z=\mp \frac{3 }{4}\lambda y \py $ and  $\mathrm{div}(Z)=\mp\frac{3}{4}\lambda$.
Since  $\deg(g_2)=2p=\rmd- 0\cdot (\rmd-\omega_1-\omega_2)$,
 by Proposition~\ref{Prop: final computation}, put $$Z_0=Z+\frac{\mathrm{div}(Z)}{\rmd-\omega_1-\omega_2- \mathrm{deg}(Z)} W=\mp \frac{3 }{2 }\lambda( x \px +   y \py) $$ and we get
   $$ g_2  \px\wedge  \py =-2Z_0(h)f\px\wedge \py=\pm  3\lambda^2 y f\px\wedge  \py= \pm  3\lambda^2 e_1f \px\wedge  \py=\pm  \frac{3}{4}\lambda^2 u\wedge v_1 \in \PH^2_{\Pi}. $$

As $p=2$,   $\deg(g_1)=p+1>2(\rmd-\omega_1-\omega_2)=2p-2$, then  $g_1=-2y^3=Y(f)\in I_f$ with $ Y=\mp \frac{1}{2}    \py.$
 Then   $\mathrm{div}(Y)=0$.
  Since $\deg(g_1)=\rmd- 1\cdot (\rmd-\omega_1-\omega_2)$,  by Proposition~\ref{Prop: final computation},
  write  $Y_0=Y+\frac{\mathrm{div}(Y)}{\rmd-\omega_1-\omega_2- \mathrm{deg}(Y)} W=Y=\mp \frac{1}{2}    \py$ and $Y_1=hY_0-\frac{2}{\deg(Y_0)} Y_0(h) W=\mp\lambda (2x\px+ \frac{3}{2} y\py).$
   Then $ g_1  \px\wedge  \py =-2Y_1(h)f\px\wedge \py=   \pm \frac{ 3}{4} \lambda^2 u\wedge v_1\in \PH^2_{\Pi}. $

 We have shown that
 if $k=1$ and $p=2$, $[v_1, w_3]_\sn=\pm \frac{ 3}{2} \lambda^2 u\wedge v_1 \in \PH^2_{\Pi}.$

\medskip

  \textbf{Case $k=2, p=2$.}  In this case,    $j=2,u_2=y$,   $g_1=-3y^2=-3u_3, g_2=-4\lambda y^{3}$. Then $ g_2=Z(f)\in I_f$ with $ Z=\mp  \lambda  \py $
   and $\mathrm{div}(Z)=0$.
 Since $\deg(g_2)=\rmd- 1\cdot (\rmd-\omega_1-\omega_2)$,
 by Proposition~\ref{Prop: final computation},
 put $Z_0=\mp  \lambda  \py$ and
  $Z_1= \mp  \lambda^2 (4  x \px +3 y \py) .$
  We have
  $$ g_2  \px\wedge  \py =-2Z_1(h)f\px\wedge \py= \pm\frac{3}{2}  \lambda^3     u\wedge v_1 \in \PH^2_{\Pi}. $$
  We obtain $[v_1, w_2]_{\sn}=-3w_3 \pm\frac{3}{2}  \lambda^3     u\wedge v_1 \in \PH^2_{\Pi}. $
\medskip

  \textbf{Case $k=3, p=2$.}  In this case,    $j=1$,   $g_1=-4y=-4u_2, g_2=-5\lambda y^{2}=-5\lambda u_3$.
  We obtain   $[v_1, w_1]_{\sn}= -4w_2-5\lambda w_3\in \PH^2_{\Pi}. $

\medskip

\textbf{Case $k=0, p=3$.} In this case, $j=5, u_j=y^{4},    g_1=-2y^{6}, g_2=-4\lambda y^{8}$.  Then $ g_1=Y(f), g_2=Z(f)\in I_f$ with  $Y=\mp \frac{1 }{3}y\py$ and $ Z=\mp \frac{2}{3}\lambda y^3\py$ and
  $\mathrm{div}(Y)=  \mp \frac{1}{3}, \mathrm{div}(Z)=  \mp 2\lambda y^{2}$.
  Let $Y_0=\mp \frac{1}{2}(x\px+y\py).$
  Since $\deg(g_1)=2p=\rmd- 0\cdot (\rmd-\omega_1-\omega_2)$ and $\deg(g_2)=3p-1=\rmd- (-1)\cdot (\rmd-\omega_1-\omega_2)$,
 by Proposition~\ref{Prop: final computation},
$$[g_1\px\wedge \py, w_{5}]_{\sn}=(-2Y_0(h)f)\px\wedge \py= \pm 2\lambda y^{2} f\px\wedge \py\in \PH^2_{\Pi}, $$
and $$[g_2\px\wedge \py, w_{5}]_{\sn}=\mathrm{div}(Z)f\px\wedge \py=\mp 2\lambda y^{2} f\px\wedge \py \in \PH^2_{\Pi}. $$
 So $[v_1, w_{5}]_{\sn}=(-2Y_0(h)f+\mathrm{div}(Z)f)\px\wedge \py=0\in \PH^2_{\Pi}. $

\medskip

 \textbf{Case $k=1, p=3$.} In this case, $j=3, u_3=y^2$,   $g_1=-4y^{4}, g_2=-6\lambda y^{6}$.

  We have  $ g_2=Z(f)\in I_f$.
 with $ Z=\mp  \lambda y \py $ and  $\mathrm{div}(Z)=\mp\lambda$.
Since  $\deg(g_2)=2p=\rmd- 0\cdot (\rmd-\omega_1-\omega_2)$,
 by Proposition~\ref{Prop: final computation}, put $Z_0=\mp \frac{3 }{2 }\lambda( x \px +   y \py) $ and
   $ g_2  \px\wedge  \py =-2Z_0(h)f\px\wedge \py=\pm  6\lambda^2 y^{2}f\px\wedge  \py= \pm  6\lambda^2 y^{2} e_1f \px\wedge  \py=\pm  \lambda^2  u\wedge v_1 \in \PH^2_{\Pi}. $
 Since $p=3$,  $\deg(g_1)=p+1\leq 2(\rmd-\omega_1-\omega_2)=2p-2$, then $g_1=-4y^{4}=-4u_{5}\in \mathcal{M}_f$.
 So $g_1  \px\wedge  \py = -4u_{5}\px\wedge  \py = -4w_{5} \in \PH^2_{\Pi}. $

 We have shown that
   if $k=1$ and $p=3$,
 $[v_1, w_3]_\sn= -4w_{5}\pm  \lambda^2   u\wedge v_1 \in \PH^2_{\Pi}.$

\medskip

  \textbf{Case $k=2, p=3$.}  In this case,    $j=1$,   $g_1=-6y^2=-6u_3, g_2=-8\lambda y^{4}=-8\lambda u_5$.
  We obtain $[v_1, w_1]_{\sn}= -6w_3-8\lambda w_5\in \PH^2_{\Pi}. $

\medskip

 \textbf{Case $k=0, p>3$.} In this case, $j=p+2, u_j=y^{p+1},    g_1=(1- p )y^{2p}, g_2=2(1-p)\lambda y^{3p-1}$.  Then $ g_1=Y(f), g_2=Z(f)\in I_f$ with  $Y=\pm \frac{1-p}{2p}y\py$ and $ Z=\pm \frac{1-p}{p}\lambda y^p\py$ and
  $\mathrm{div}(Y)=  \pm \frac{1-p}{2p}, \mathrm{div}(Z)=  \pm \lambda(1-p)y^{p-1}$.
  Let $Y_0=\mp \frac{1}{2}(x\px+y\py).$
  Since $\deg(g_1)=2p=\rmd- 0\cdot (\rmd-\omega_1-\omega_2)$ and $\deg(g_2)=3p-1=\rmd- (-1)\cdot (\rmd-\omega_1-\omega_2)$,
 by Proposition~\ref{Prop: final computation},
$$[g_1\px\wedge \py, w_{p+2}]_{\sn}=(-2Y_0(h)f)\px\wedge \py= \pm \lambda (p-1)y^{p-1} f\px\wedge \py\in \PH^2_{\Pi}, $$
and $$[g_2\px\wedge \py, w_{p+2}]_{\sn}=\mathrm{div}(Z)f\px\wedge \py=\pm \lambda (1-p)y^{p-1} f\px\wedge \py \in \PH^2_{\Pi}. $$
 So $$[v_1, w_{p+2}]_{\sn}=(-2Y_0(h)f+\mathrm{div}(Z)f)\px\wedge \py=0\in \PH^2_{\Pi}. $$

\medskip

 \textbf{Case $k=1, p>3$.} In this case, $j=3, u_j=y^2$,   $g_1=2(1-p)y^{p+1}, g_2=3(1-p)\lambda y^{2p}$.

  We have  $ g_2=Z(f)\in I_f$.
 with $ Z=\pm \frac{3(1-p)}{2p}\lambda y \py $ and  $\mathrm{div}(Z)=\pm \frac{3(1-p)}{2p}\lambda$.
Since  $\deg(g_2)=2p=\rmd- 0\cdot (\rmd-\omega_1-\omega_2)$,
 by Proposition~\ref{Prop: final computation}, put $Z_0=\mp \frac{3 }{2 }\lambda( x \px +   y \py) $ and
   $$ g_2  \px\wedge  \py =-2Z_0(h)f\px\wedge \py=\pm  3\lambda^2(p-1)y^{p-1}f\px\wedge  \py= \pm  3\lambda^2(p-1)e_1f \px\wedge  \py=\pm 3\lambda^2\frac{p-1}{2p} u\wedge v_1 \in \PH^2_{\Pi}. $$
As $p>3$,  $\deg(g_1)=p+1< 2(\rmd-\omega_1-\omega_2)=2p-2$, then  $g_1=2(1-p)y^{p+1}=2(1-p)u_{p+2}\in \mathcal{M}_f$.
 So $$g_1  \px\wedge  \py = 2(1-p)u_{p+2}\px\wedge  \py = 2(1-p)w_{p+2} \in \PH^2_{\Pi}. $$

 We have shown that
  if $k=1$ and $p> 3$,
 $[v_1, w_3]_\sn= 2(1-p)w_{p+2}\pm 3\lambda^2\frac{p-1}{2p}  u\wedge v_1 \in \PH^2_{\Pi}.$

\medskip

 \begin{Prop}  Let $\Pi=(x^2\pm y^{2p})(1+\lambda y^{p-1})\frac{\partial}{\partial x}\wedge \frac{\partial}{\partial y}$ with $p\geq 1$ be of type $A_{2p-1}^\pm$.
\begin{itemize}
\item[(a)] If $p=1$, there exists an isomorphism of Gerstenhaber  algebras
 $\PH^{*}_{\Pi}\cong \wedge_\k(u, v_1')\times_{\mathbb{K}}
 \mathbb{K}[w_1]/(w_{1}^{2}) $
with $|u|=|v_1'|=1, |w_{1}|=2$  and the only nonzero bracket on generators is given by
 $
  [v_1', w_{1}]_{\sn}= w_1. $

 \item[(b)] If $p=2$ and $\lambda=0$, then $\Pi=\Pi_0$ and there exists an isomorphism of Gerstenhaber  algebras $$\PH^{*}_{\Pi_0}\cong\wedge_\k(u, v_{1})\times_{\k}\k[w_{1}]/(w_{1}^{2})\times_{\k}\k[w_{2}']/{(w_{2}'}^{2}) \times_{\k}\k[w_{3}']/({w_{3}'}^{2}) $$
 with $|u|=|v_1|=1, |w_{1}|=|w_2'|=|w_3'|=2$  and the nontrivial  brackets are given as follows:
  $$[v_1, w_1]_{\sn}=w_2', [v_1, w_2']_{\sn}=w_3'.$$

 If $p=2$ and $\lambda\neq 0$, there exists an isomorphism of Gerstenhaber  algebras $$\PH^{*}_{\Pi}\cong\wedge_\k(u', v_{1})\times_{\k}\k[w_{1}]/(w_{1}^{2}) \times_{\k}\k[w_{2}']/{(w_{2}'}^{2}) \times_{\k}\k[w_{3}']/({w_{3}'}^{2}) $$  with $|u'|=|v_1|=1, |w_{1}|=|w_2'|=|w_3'|=2$  and the nontrivial  brackets are given by:
 $$[v_1, w_1]_{\sn}=w_2', [v_1, w_2']_{\sn}=w_3', [v_1, w_3']=u'v_1.$$

 \item[(c)] If $p=3$ and $\lambda=0$, then $\Pi=\Pi_0$ and  there exists an isomorphism of Gerstenhaber  algebras $$\PH^{*}_{\Pi_0}\cong\wedge_\k(u, v_{1})\times_{\k}\k[w_{1}]/(w_{1}^{2})\times_{\k} \k[w_{2}]/(w_{2}^{2})\times_{\k}\k[w_{3}']/({w_{3}'}^{2})\times_{\k}\k[w_{4}']/({w_{4}'}^{2})\times_{\k}\k[w_{5}']/({w_{5}'}^{2}) $$
     with $|u|=|v_1|=1, |w_{1}|=|w_{2}|=|w_{3}'|=|w_{4}'| =|w_5'|=2$  and the nontrivial brackets are given by
    $$[v_1, w_1]_{\sn}=w_3', [v_1, w_3']_{\sn}=w_5', [v_1, w_2]_{\sn}=w_4'.$$

    If $p=3$ and $\lambda\neq 0$, then   there exists an isomorphism of Gerstenhaber  algebras $$\PH^{*}_{\Pi}\cong\wedge_\k(u, v_{1})\times_{\k}\k[w_{1}]/(w_{1}^{2})\times_{\k} \k[w_{2}]/(w_{2}^{2})\times_{\k}\k[w_{3}']/({w_{3}'}^{2})\times_{\k}\k[w_{4}']/({w_{4}'}^{2})\times_{\k}\k[w_{5}']/({w_{5}'}^{2}) $$
     with $|u|=|v_1|=1, |w_{1}|=|w_{2}|=|w_{3}'|=|w_{4}'| =|w_5'|=2$  and the nontrivial brackets are given by
    $$[v_1, w_1]_{\sn}=w_3', [v_1, w_3']_{\sn}=w_5', [v_1, w_2]_{\sn}=w_4'.$$

  \item[(d)]If $p\geq 4$ and $\lambda=0$, then $\Pi=\Pi_0$ and there exists an isomorphism of Gerstenhaber  algebras
     $$\PH^{*}_{\Pi_0}\cong\wedge_\k(u, v_{1})\times_{\k}\k[w_{1}]/(w_{1}^{2})\times_{\k} \cdots \times_{\k}\k[w_{p-1}]/(w_{p-1}^{2})
     \times_{\k}\k[w_{p}']/({w_{p}'}^{2})\times_{\k}  \cdots \times_{\k}\k[w_{2p-1}']/({w_{2p-1}'}^{2}) $$
  with $|u|=|v_1|=1, |w_{1}|=\cdots=|w_{p-1}|=|w_p'|=\cdots=|w_{2p-1}'|=2$  and the nontrivial brackets are given by
  $$[v_1, w_1]_\sn=w_p', [v_1, w_p']_\sn=w_{2p-1}', [v_1, w_j]_{\sn}=w_{p+j-1}',j=2, \cdots, p-1.$$

  If $p\geq 4$ and $\lambda\neq 0$, then   there exists an isomorphism of Gerstenhaber  algebras
     $$\PH^{*}_{\Pi}\cong\wedge_\k(u, v_{1})\times_{\k}\k[w_{1}]/(w_{1}^{2})\times_{\k} \cdots \times_{\k}\k[w_{p-1}]/(w_{p-1}^{2})
     \times_{\k}\k[w_{p}']/({w_{p}'}^{2})\times_{\k}  \cdots \times_{\k}\k[w_{2p-1}']/({w_{2p-1}'}^{2}) $$
  with $|u|=|v_1|=1, |w_{1}|=\cdots=|w_{p-1}|=|w_p'|=\cdots=|w_{2p-1}'|=2$  and the nontrivial brackets are given by
  $$[v_1, w_1]_\sn=w_p', [v_1, w_p']_\sn=w_{2p-1}', [v_1, w_j]_{\sn}=w_{p+j-1}',j=2, \cdots, p-1.$$
 \end{itemize}

 \end{Prop}

 \begin{Proof}

 (a) It suffices to replace $v_1$ by $v_1':=-\frac{1}{2}v_1$.

 (b) By the above discussion, we obtain that
  $$ [v_1, w_1]_\sn  = -4w_2-5\lambda w_3,
  [v_1, w_2]_\sn  = -3w_3  \pm\frac{3}{2}\lambda^3 u \wedge  v_1,
   [v_1, w_3]_\sn  = \pm\frac{3}{2}\lambda^2   u\wedge v_1,
  [v_1, w_4]_\sn = 0\in \PH^2_{\Pi}.$$
Define $
   w_2'= -4w_2-5\lambda w_3, w_3'=[v_1, w_2']_\sn  = 12w_3 \mp \frac{27}{2}\lambda^3  u\wedge v_1.$
  Then $[v_1, w_3']_\sn= \pm  18 \lambda^2 u\wedge v_1$.
   If $\lambda\neq 0$, put $u'=\pm  18 \lambda^2 u,$ otherwise, do not change $u$,
   then we have the desired result.

(c) By the above discussion, we obtain that
  $ [v_1, w_1]_{\sn}  =  -6w_3-8\lambda w_5,
  [v_1, w_3]_{\sn}  =  -4w_5   \pm \lambda^2  u \wedge  v_1, [v_1, w_5]_\sn=  0\in \PH^2_{\Pi}.
   $

   It remains to compute $[v_1, w_2]_\sn$ and $[v_1, w_4]_\sn$.
 For $j=2$, $g_1=-5y^3=-5u_4$ and   $ g_1  \px\wedge  \py =-5 w_4 \in \PH^2_{\Pi}. $
   $g_2=-7\lambda y^5=Z(f)\in I_f$, so $ g_2  \px\wedge  \py =0 \in \PH^2_{\Pi}. $
   We obtain  $[v_1, w_2]_\sn=-5w_4\in \PH^2_{\Pi}.$
For $j=4$, $g_1=-3y^6, g_2=-5\lambda y^7\in I_f$, so  $[v_1, w_4]_\sn=0\in \PH^2_{\Pi}$.

  Put
  $w_3'=  -6w_3-8\lambda w_5, w_4'=-5w_4, w_5'= 24w_5   \mp 6\lambda^2  u \wedge  v_1, $ then   the   result follows.

  (d) By the above discussion, we obtain that $ [v_1, w_{p+2}]_{\sn}  =0$ and that
  $$ [v_1, w_3]_{\sn}  =  2(1-p)w_{p+2}\pm 3\lambda^2\frac{p-1}{2p}  u\wedge v_1 \in \PH^2_{\Pi}.$$
  It remains to compute $[v_1, w_j]_\sn$ for $j\neq 3, p+2$.

  For $j=1$, $g_1=-2py^{p-1}=-2pu_{p}, g_2=(1-3p)\lambda y^{2p-2} =(1-3p)\lambda u_{2p-1}$. So by Proposition~\ref{Prop: final computation}, $$[v_1, w_1]_\sn=-2pw_{p}+(1-3p)\lambda w_{2p-1}\in \PH^2_{\Pi}.$$

  For $j=2$ or $4\leq j\leq p$, $g_1=(j-2p-1)y^{p+j-2}=(j-2p-1)u_{p+j-1}, g_2=(j-3p)\lambda y^{2p+j-3}\in I_f$. So by Proposition~\ref{Prop: final computation}, $$[v_1, w_j]_\sn=(j-2p-1)w_{p+j-1}\in \PH^2_{\Pi}.$$

  For $j=p+1$ or $p+3\leq j\leq 2p-1$, $g_1=(j-2p-1)y^{p+j-2}, g_2=(j-3p)\lambda y^{2p+j-3}\in I_f$. So by Proposition~\ref{Prop: final computation}, $$[v_1, w_j]_\sn=0\in \PH^2_{\Pi}.$$

 Put $w_p'=-2pw_{p}+(1-3p)\lambda w_{2p-1}, w_{p+1}'=(1-2p)w_{p+1}, w_{p+2}'= 2(1-p)w_{p+2}\pm  \frac{3(p-1)}{2p}   \lambda^2 u\wedge v_1$ and $ w_j'=(j-3p)w_{j}, j=p+3, \cdots, 2p-2 $ and $ w_{2p-1}'=2p(2p+1)w_{2p-1}$,   then we obtain the desired result.
 \end{Proof}





\subsection{Type  $D_{2p}^{\pm}$: $(x^{2}y\pm y^{2p-1})(1+\lambda x+\mu y^{p-1}),  p\ge 2$}

In this case, $f=x^{2}y\pm y^{2p-1}, h=\lambda x+\mu y^{p-1}$, $\omega_{1}=p-1,~\omega_{2}=1, d=2p-1$ and $
\rmd- \omega_{1}-\omega_{2}=p-1$. A basis of $\mathcal{P}_{\rmd- \omega_{1}-\omega_{2}}=\mathcal{P}_{p-1}$ is given by  $\{e_1:=y^{p-1}, e_2:=x\}$ and $r=2$. Since $I_{f}=(\frac{\partial f}{\partial x},\frac{\partial f}{\partial y})=( 2xy,x^{2}\pm (2p-1)y^{2p-2})$,   then  $c=2p$ and a monomial basis of $\mathcal{M}_f=\mathcal{F} /I_{f}$ can be chosen as  $$ u_1:=1,   u_2:=y, u_3:=y^{2},\cdots, u_{2p-2}:=y^{2p-3}, u_{2p-1}:=y^{2p-2}, u_{2p}=x$$
whose degrees are respectively $0, 1, \cdots, 2p-2, p-1$.

    Let  $u=(1+h)H_{f}=(1+\lambda x+\mu y^{p-1})((x^{2}\pm (2p-1)y^{2p-2})\frac{\partial}{\partial x}-2xy\frac{\partial}{\partial y}),$
    $v_1=(1+h)e_{1}W=(1+\lambda x+\mu y^{p-1})y^{p-1}((p-1)x\frac{\partial}{\partial x}+y\frac{\partial}{\partial y})$ and
    $v_2=(1+h)e_{2}W=(1+\lambda x+\mu y^{p-1})x((p-1)x\frac{\partial}{\partial x}+y\frac{\partial}{\partial y}).$
A basis of $\mathcal{M}_f\frac{\partial }{\partial x}\wedge \frac{\partial }{\partial y} $ is given by $u_{1}\frac{\partial}{\partial x}\wedge\frac{\partial}{\partial y},\cdots, u_{2p}\frac{\partial}{\partial x}\wedge\frac{\partial}{\partial y}$,  and $e_{1}f\frac{\partial}{\partial x}\wedge\frac{\partial}{\partial y} , e_{2}f\frac{\partial}{\partial x}\wedge\frac{\partial}{\partial y}$  form a basis of $\mathcal{P}_{\rmd-\omega_1-\omega_2}f \frac{\partial }{\partial x}\wedge \frac{\partial }{\partial y}$. Denote
$w_j=u_{j}\frac{\partial}{\partial x}\wedge\frac{\partial}{\partial y}, j=1, \cdots, 2p$.

 By Theorem~\ref{Thm: wedge product result}, we have an isomorphism of graded algebras
 $$\PH^{*}_{\Pi}\cong \k\langle u, v_1,v_{2}\rangle/(u^{2}, v_1^2,   v_{2}^{2},u v_{1}+v_{1}u, u v_{2}+v_{2}u, v_{1}v_{2})\times_{\k}\k[w_{1}]/(w_{1}^{2})\times_{\k}\cdots\times_{\k}\k[w_{2p}]/(w_{2p}^{2}).$$

By Subsection~\ref{subsection bracket for Pi}, to compute the bracket, we need to consider
  $[v_i, w_j]_{\sn}, i=1, 2,  j=1, \cdots, 2p$.

\subsubsection{$[v_1, w_j]$}

  By Lemma~\ref{Lem: [v_i, w_j]},  for $1\leq j\leq 2p$,
$$\begin{array}{rcl}
[v_1, w_j]_{\sn}&=&[(1+h)e_{1}W,u_{j}\frac{\partial}{\partial x}\wedge\frac{\partial}{\partial y}]_{\sn}\\
&=& ((\deg(u_j)-2p+1)y^{p-1}u_j+( \deg(u_j)-3p+2)(\lambda x+\mu y^{p-1})y^{p-1}u_j )\frac{\partial}{\partial x}\wedge \frac{\partial}{\partial y}
\end{array}$$
 As in Lemma~\ref{Lem: [v_i, w_j]}, let $g_1=(\deg(u_j)-2p+1)y^{p-1}u_j$ and $g_2=( \deg(u_j)-3p+2)(\lambda x+\mu y^{p-1})y^{p-1}u_j$.

 It is easy to see that  $\deg(g_1)=p-1+\deg(u_j)=\rmd- k\cdot (\rmd-\omega_1-\omega_2)=2p-1-k(p-1)$ with $k\geq -1$ if and only if
 $k=0, 1$ or ($k=2$ and $p=2$).

\medskip

 Let $k=0$. Then $\deg(u_j)=p, u_j=y^p, j=p+1, g_1=(1-p)y^{2p-1}.$ Since $\deg(g_1)=2p-1>2(\rmd-\omega_1-\omega_2)=
2p-2$, $g_1\in I_f$ and $g_1=Y(f)$ with
$Y=\pm \frac{p-1}{2(2p-1)}(x\px-2y\py)$.
 Let $Y_0= \mp \frac{1}{2}y\py.$
 By Proposition~\ref{Prop: final computation}, we obtain that
 $$ g_1\px\wedge \py =-2Y_0(h)f\px\wedge \py=\pm \mu (p-1)y^{p-1}f\px\wedge \py=\pm \frac{p-1}{2p-1}\mu u\wedge v_1\in \PH^2_\Pi.$$

Let $k=1$.  Then $\deg(u_j)=1$.   If $p>2$, then $  u_j=y , j=2, g_1=2(1-p)y^{p}=2(1-p)u_{p+1}.$
 By Proposition~\ref{Prop: final computation}, we see that
 $$ g_1\px\wedge \py = 2(1-p)w_{p+1}\in \PH^2_\Pi.$$
 If $p=2$, then either  $  u_j=y , j=2, g_1=-2y^2=-2u_{3},$ or  $  u_j=x , j=4, g_1=-2xy.$
In the first case,  by Proposition~\ref{Prop: final computation}, we see
 $$ g_1\px\wedge \py = -2u_3\in \PH^2_\Pi;$$
 in the second case,   $g_1=Y(f)$ with
$Y=-\px$.
 Let  $Y_0=-\px $  and
   $Y_1= -(3\lambda x+\mu y)\px- 2\lambda y\py.$
   By Proposition~\ref{Prop: final computation},   $$  g_1  \px\wedge  \py  =-2Y_1(h)f\px\wedge \py= 6\lambda (\lambda x+\mu y)f\px\wedge \py=   2\lambda (\mu u\wedge v_1
    + \lambda  u \wedge v_2)\in \PH^2_{\Pi}. $$

 Let $k=2$. Then $p=2, \deg(u_j)=0, u_j=1 , j=1, g_1= -3y = -3 u_{2}.$   By Proposition~\ref{Prop: final computation}, we see that
 $$  g_1\px\wedge \py = -3w_{2}\in \PH^2_\Pi.$$

Suppose that there is no $k\geq -1$ such that $\deg(g_1)=p-1+\deg(u_j)=\rmd- k\cdot (\rmd-\omega_1-\omega_2)=2p-1-k(p-1)$.
Then   $p>2$ and $j\neq p+1, 2$.
If $j=1$ or $3\leq j\leq p$, then $g_1=(j-2p )y^{p+j-2}=(j-2p)u_{p+j-1} $ and by Proposition~\ref{Prop: final computation},
$$ g_1\px\wedge \py = (j-2p)w_{p+j-1}\in \PH^2_\Pi;$$
if $p+2\leq j\leq 2p-1$, then $g_1=(j-2p)y^{p+j-2}\in I_f$,  by Proposition~\ref{Prop: final computation},
$$ g_1\px\wedge \py = 0\in \PH^2_\Pi,$$
because $p+j-2>2(\rmd-\omega_1-\omega_2)=2p-2$ and $g_1\in I_f$; if $j=2p$, then $g_1=-pxy^{p-1}\in I_f$ and   by Proposition~\ref{Prop: final computation}, $$ g_1\px\wedge \py = 0\in \PH^2_\Pi.$$

\medskip

We have proved the following:
If $p=2$,
$$ g_1\px\wedge \py = \left\{\begin{array}{ll} -3w_2, &j =1, \\  -2w_3,&  j=2, \\  \pm \frac{ 1}{3 }\mu u\wedge v_1,  & j=3, \\  2\lambda (\mu u\wedge v_1
    + \lambda  u \wedge v_2),  & j=4;\end{array}\right.$$
if $p>2$,
$$ g_1\px\wedge \py = \left\{\begin{array}{ll}  (1-2p)w_p, &j =1, \\   2(1-p)w_{p+1},& j =2, \\   (j-2p)w_{p+j-1},  & 3\leq j\leq p, \\  \pm \frac{p-1}{2p-1}\mu u\wedge v_1,  & j=p+1,\\
0& p+1\leq j\leq 2p.\end{array}\right.$$

\medskip

Now consider $g_2=( \deg(u_j)-3p+2)(\lambda x+\mu y^{p-1})y^{p-1}u_j$.
It is easy to see that  $\deg(g_2)=2(p-1)+\deg(u_j)=\rmd- \ell \cdot (\rmd-\omega_1-\omega_2)=2p-1-\ell (p-1)$ with $\ell\geq -1$ if and only if
 $\ell=-1, 0, 1$ and $\ell=1$ only when $p=2$.

Let $\ell=1$.  Then $p=2$ and $\deg(u_j)=0, u_j=1, j=1$. We see $g_2=-4(\lambda x+\mu y)y=-4\mu u_3+Z(f)$ with $Z=-2\lambda \px$. Let  $Z_0=-2\lambda \px $
and
   $Z_1=  (-6\lambda^2x-2\lambda\mu y ) \px-4\lambda^2 y\py.$
   Then $$  g_2  \px\wedge  \py  =-4\mu w_3-2Z_1(h)f\px\wedge \py=-4\mu w_3 +4 \lambda^2\mu u\wedge v_1+4 \lambda^3 u\wedge v_2  \in \PH^2_{\Pi}. $$

 Let $\ell=0$. Then $\deg(u_j)=1$. There are two cases: either  $u_j=y , j=2, g_2=3(1-p)(\lambda x+\mu y^{p-1})y^p, $
 or $p=2, j=4, u_j=x, g_2=3(1-p)(\lambda x+\mu y)xy$.

 In the first case,
 since $\deg(g_2)=2p-1>2(\rmd-\omega_1-\omega_2)=
2p-2$, $g_2\in I_f$ and $g_2=Z(f)$ with
$$Z= (\frac{3}{2}(1-p)\lambda y^{p-1}\pm \frac{3(p-1)}{2(2p-1)}\mu x)\px  \mp \frac{3(p-1)}{2p-1}\mu  y\py .$$
 Let $Z_0=\frac{3}{2}(1-p)\lambda y^{p-1}\px\mp \frac{3 }{2 }\mu y\py  .$
 We obtain that
 $$ g_2\px\wedge \py =-2Z_0(h)f\px\wedge \py =  \frac{3(p-1)}{2p-1}(\lambda^2 \pm     \mu^2)u\wedge v_1 \in \PH^2_\Pi.$$

  In the second case,  $p=2, j=4, u_j=x, g_2=-3(\lambda x+\mu y)xy$. Since $\deg(g_2)=3>2(\rmd-\omega_1-\omega_2)=
2$, $g_2\in I_f$ and $g_2=Z(f)$ with $Z=-\frac{3}{2} (\lambda x+\mu y)\px$.
Let $Z_0=(-3 \lambda x-\frac{3}{2} \mu y)\px -\frac{3}{2 }\lambda y\py  .$
 We obtain that
 $$[g_2\px\wedge \py, w_{4}]_{\sn}=-2Z_0(h)f\px\wedge \py= (6\lambda^2 x+ 6\lambda \mu y)f\px\wedge \py
  =2\lambda \mu u\wedge v_1+ 2\lambda^2 u\wedge v_2     \in \PH^2_\Pi.$$

Let $\ell=-1$. Then $\deg(u_j)=p, u_j=y^p , j=p+1, g_2= 2(1-p)(\lambda x+\mu y^{p-1})y^{2p-1}.$
 Since $\deg(g_2)=3p-2>2(\rmd-\omega_1-\omega_2)=
2p-2$, $g_2\in I_f$ and $g_2=Z(f)$ with
$Z= (1-p)  y^{p-1} (\lambda y^{p-1}\mp \frac{1}{2p-1}\mu x)\px \pm \frac{2(1-p)}{2p-1}\mu  y^p\py .$
  We obtain that
 $$[g_2\px\wedge \py, w_{p+1}]_{\sn}=\mathrm{div}(Z)f\px\wedge \py= \pm (1-p)\mu e_1f\px\wedge \py=
 \pm \frac{1-p}{ 2p-1  }\mu u\wedge v_1 \in \PH^2_\Pi.$$

Suppose that there is no $\ell\geq -1$ such that $g_2=( \deg(u_j)-3p+2)(\lambda x+\mu y^{p-1})y^{p-1}u_j =\rmd- k\cdot (\rmd-\omega_1-\omega_2)=2p-1-k(p-1)$.
Then   $p>2$ and $j\neq p+1, 2$.
As  $g_2=( \deg(u_j)-3p+2)(\lambda x+\mu y^{p-1})y^{p-1}u_j$,
If $j=1$, then $g_2=(2-3p)(\lambda x+\mu y^{p-1})y^{p-1}=(2-3p)\mu u_{2p-1}+Z(f)$ with $Z$ a vector field, and
$$ g_2\px\wedge \py = (2-3p)\mu w_{2p-1} \in \PH^2_\Pi;$$
if $3\leq j\leq p$, then $g_2=(j-3p+1 )(\lambda x+\mu y^{p-1})y^{p+j-2}= Z(f)$ for a vector field $Z$ and
$$ g_2\px\wedge \py =0\in \PH^2_\Pi;$$
if $p+2\leq j\leq 2p-1$, then $g_2=(j-3p+1 )(\lambda x+\mu y^{p-1})y^{p+j-2}
= Z(f)$ for a vector field $Z$ and
$$ g_2\px\wedge \py =0\in \PH^2_\Pi;$$
if $j=2p$, $g_2=(1-2p )(\lambda x+\mu y^{p-1})xy^{p-1}
= Z(f)$ for a vector field $Z$  and
$$ g_2\px\wedge \py =0\in \PH^2_\Pi.$$

\medskip

We have proved the following:
If $p=2$,

$$ g_2\px\wedge \py = \left\{\begin{array}{ll} -4\mu w_3 +4 \lambda^2\mu u\wedge v_1+4 \lambda^3 u\wedge v_2 , &j =1, \\
    (\lambda^2 \pm     \mu^2)u\wedge v_1,&  j=2, \\
       \mp \frac{1 }{ 3  }\mu u\wedge v_1,  & j=3, \\
       2\lambda \mu u\wedge v_1+ 2\lambda^2 u\wedge v_2,  & j=4.\end{array}\right.$$

If $p>2$,
$$ g_2\px\wedge \py = \left\{\begin{array}{ll}  (2-3p)\mu w_{2p-1}, &j =1, \\
   \frac{3(p-1)}{2p-1}(\lambda^2 \pm     \mu^2)u\wedge v_1,& j =2, \\
      0,  & 3\leq j\leq p, \\
         \pm \frac{1-p}{ 2p-1  }\mu u\wedge v_1,  & j=p+1,\\
0& p+1\leq j\leq 2p.\end{array}\right.$$

\medskip

 Hence, we obtain $[v_1, w_j]_\sn$:
If $p=2$,
$$[v_1, w_{j}]_{\sn}= \left\{\begin{array}{ll}  -3w_2-4\mu w_3 +4 \lambda^2\mu u\wedge v_1+4 \lambda^3 u\wedge v_2 , &j =1, \\
    -2w_3+(\lambda^2 \pm     \mu^2)u\wedge v_1,&  j=2, \\
      0,  & j=3, \\
       4\lambda \mu u\wedge v_1+ 4\lambda^2 u\wedge v_2,  & j=4.\end{array}\right.$$
if $p>2$,
$$[v_1, w_{j}]_{\sn}= \left\{\begin{array}{ll} (1-2p)w_p+ (2-3p)\mu w_{2p-1}, &j =1, \\
  2(1-p)w_{p+1}+ \frac{3(p-1)}{2p-1}(\lambda^2 \pm     \mu^2)u\wedge v_1,& j =2, \\
      (j-2p)w_{p+j-1},  & 3\leq j\leq p, \\

0& p+1\leq j\leq 2p.\end{array}\right.$$

\medskip

\medskip
\subsubsection{$[v_2, w_j]$}

 By Lemma~\ref{Lem: [v_i, w_j]},  for $1\leq j\leq 2p-1$,
$$\begin{array}{rcl}
[v_2, w_j]_{\sn}&=&[(1+h)e_{2}W,u_{j}\frac{\partial}{\partial x}\wedge\frac{\partial}{\partial y}]_{\sn}\\
&=& ((\deg(u_j)-2p+1)xu_j+( \deg(u_j)-3p+2)(\lambda x+\mu y^{p-1})xu_j )\frac{\partial}{\partial x}\wedge \frac{\partial}{\partial y}
\end{array}$$
 As in Lemma~\ref{Lem: [v_i, w_j]}, let $g_1'=(\deg(u_j)-2p+1)xu_j$ and $g_2'=( \deg(u_j)-3p+2)(\lambda x+\mu y^{p-1})xu_j$.

 It is easy to see that  $\deg(g_1')=p-1+\deg(u_j)=\rmd- k\cdot (\rmd-\omega_1-\omega_2)=2p-1-k(p-1)$ with $k\geq -1$ if and only if
 $k=0, 1,$ or $2$ (and $p=2)$.

\medskip

 Let $k=0$. Then $\deg(u_j)=p, u_j=y^p, j=p+1, g_1'=(1-p)xy^{p}.$ Since $\deg(g_1')=2p-1>2(\rmd-\omega_1-\omega_2)=
2p-2$, $g_1'\in I_f$ and $g_1'=Y(f)$ with
$Y=  \frac{1-p}{2}y^{p-1}\px$.
 Let $Y_0= \frac{1-p}{2}y^{p-1}\px.$
 We obtain that
 $$ g_1'\px\wedge \py =-2Y_0(h)f\px\wedge \py=  (p-1)\lambda y^{p-1}f\px\wedge \py=
 \frac{p-1}{2p-1} \lambda u\wedge v_1\in \PH^2_\Pi.$$

Let $k=1$.  Then $\deg(u_j)=1$.

  If $p>2$, then $  u_j=y , j=2, g_1'=2(1-p) xy=Y(f)$ with $Y=(1-p)\px$.
  Let $Y_0=(1-p)\px $
  and $Y_1=(3(1-p)\lambda  x+(1-p)\mu y^{p-1} )\px-2\lambda y\py.$
 We obtain that
 $$ g_1'\px\wedge \py =-2Y_1(h)f\px\wedge \py=     \frac{6(p-1) }{2p-1}\lambda \mu u\wedge v_1+\frac{6(p-1) }{2p-1}\lambda^2 u\wedge v_2\in \PH^2_\Pi.$$

 If $p=2$, then either  $  u_j=y , j=2, g_1'=-2xy=Y(f) $ with $Y=-\px$, or  $  u_j=x , j=4,
 g_1'=-2x^2=\pm 6y^2-2(x^2\pm 3y^2)=\pm 6u_{3}+Y(f)$ with $Y=-2\py$.

In the first case, let $Y_0=- \px $
  and $Y_1=( -3\lambda   x-\mu y )\px- 2\lambda y\py.$
 We obtain that
 $$ g_1'\px\wedge \py = -2Y_1(h)f\px\wedge \py =2\lambda \mu   u\wedge v_1+  2\lambda^2 u\wedge v_2 \in \PH^2_\Pi.$$

 In the second case,   let  $Y_0=-2 \py $  and
   $Y_1=  (-4 \mu x)\px+(-2\lambda x-6\mu y)\py.$
   Then $$    g_1'  \px\wedge  \py  =\pm 6w_{3}-2Y_1(h)f\px\wedge \py= \pm 6w_{3} +4\mu^2 u\wedge v_1+4\lambda \mu  u\wedge v_2 \in \PH^2_{\Pi}. $$

 Let $k=2$. Then $p=2, \deg(u_j)=0, u_j=1 , j=1, g_1'= -3x = -3 u_{4}.$   We see that
 $$ g_1'\px\wedge \py = -3w_{4}\in \PH^2_\Pi.$$

Suppose that there is no $k\geq -1$ such that $\deg(g_1')=p-1+\deg(u_j)=\rmd- k\cdot (\rmd-\omega_1-\omega_2)=2p-1-k(p-1)$.
Then   $p>2$ and $j\neq p+1, 2$.

 If $j=1$, then $g_1'=(1-2p)x=(1-2p)u_{2p}$ and   $$ g_1'\px\wedge \py = (1-2p)w_{2p}\in \PH^2_\Pi.$$

If $3\leq j\leq p$, then $g_1'=(j-2p)xy^{j-1}\in I_f$ and
$$ g_1'\px\wedge \py = 0\in \PH^2_\Pi.$$

If $p+2\leq j\leq 2p-1$, then $g_1'=(j-2p)xy^{j-1}\in I_f$,
$$ g_1'\px\wedge \py = 0\in \PH^2_\Pi.$$

If $j=2p$, then $g_1'=  -p x^2=\pm  p(2p-1)u_{2p-1}\, \mathrm{mod}\, I_f$ and $$ g_1'\px\wedge \py =
 \pm p(2p-1) w_{2p-1}\in \PH^2_\Pi.$$

\medskip

We have proved the following:
If $p=2$,
$$ g_1'\px\wedge \py = \left\{\begin{array}{ll} -3w_{4}, &j =1, \\    2  \lambda \mu   u\wedge v_1+  2\lambda^2 u\wedge v_2,&  j=2, \\    \frac{1}{3} \lambda u\wedge v_1,  & j=3, \\ \pm 6w_{3} +4\mu^2 u\wedge v_1+4\lambda \mu  u\wedge v_2,  & j=4;\end{array}\right.$$
if $p>2$,
$$ g_1'\px\wedge \py = \left\{\begin{array}{ll}  (1-2p)w_{2p}, &j =1, \\
 \frac{6(p-1)}{2p-1} \lambda \mu u\wedge v_1+\frac{6(p-1)}{2p-1} \lambda^2 u\wedge v_2,& j =2, \\   0,  & 3\leq j\leq p, \\    \frac{p-1}{2p-1} \lambda u\wedge v_1,  & j=p+1,\\
0& p+1\leq j\leq 2p-1\\
 \pm p(2p-1) w_{2p-1},& j=2p.\end{array}\right.$$

\medskip

Now consider $g_2'=( \deg(u_j)-3p+2)(\lambda x+\mu y^{p-1})xu_j$.
It is easy to see that  $\deg(g_2')=2(p-1)+\deg(u_j)=\rmd- \ell \cdot (\rmd-\omega_1-\omega_2)=2p-1-\ell (p-1)$ with $\ell \geq -1$ if and only if
 $\ell=-1, 0, 1$ and $\ell=1$ only when $p=2$.

Let $\ell=1$.  Then $p=2$ and $\deg(u_j)=0, u_j=1, j=1$. We see $g_2'=-4\lambda x^2-4\mu xy =\pm 12\lambda  u_{3}+Z(f)  $ with $Z=-2\mu \px-4\lambda \py$.
Let  $Z_0=-2\mu \px-4\lambda \py$
and
   $Z_1= ( -14\lambda\mu x-2\mu^2 y)\px +(-4\lambda^2x -16\lambda \mu y)     \py.$
   Then $$ \begin{array}{rcl}
   g_2'  \px\wedge  \py  &=&\pm 12\lambda w_{2p-1}-2Z_1(h)f\px\wedge \py\\

    &=&\pm 12\lambda w_{3} + 12\lambda \mu^2  u\wedge v_1+12\lambda^2 \mu\wedge v_2  \in \PH^2_{\Pi}.\end{array} $$

 Let $\ell=0$. Then $\deg(u_j)=1$. There are two cases: either  $u_j=y , j=2, g_2'=3(1-p)(\lambda x+\mu y^{p-1})xy, $
 or $p=2, j=4, u_j=x, g_2'=-3(\lambda x+\mu y)x^2$.

 In the first case,
 since $\deg(g_2')=2p-1>2(\rmd-\omega_1-\omega_2)=
2p-2$, $g_2'\in I_f$ and $g_2'=Z(f)$ with
$Z=  \frac{3}{2}(1-p)(\lambda x+\mu y^{p-1})\px$.
 Let $Z_0=(  3(1-p)\lambda  x+ \frac{3}{2}(1-p)\mu y^{p-1})\px  - \frac{3}{2}\lambda y\py.$
 We obtain that $$\begin{array}{rcl}
   g_2'\px\wedge \py &=&-2Z_0(h)f\px\wedge \py\\
  &=& (6(p-1)\lambda^2 x + 6\lambda \mu(p-1)   y^{p-1})f\px\wedge \py\\
  &=&   \frac{6(p-1)}{2p-1}\lambda \mu      u\wedge v_1+  \frac{6(p-1)}{2p-1}\lambda^2  u\wedge v_2 \in \PH^2_\Pi.
  \end{array}$$

  In the second case,  $p=2, j=4, u_j=x$ and $  g_2'=-3(\lambda x+\mu y)x^2= Z(f)$ with $Z=(-\frac{3}{2}
\mu  x\pm \frac{9}{2} \lambda  y)\px -3\lambda x\py.$
Let $Z_0= (-3\mu x\pm  \frac{9}{2}\lambda y)\px- (3\lambda x+\frac{3}{3}\mu y)\py.$
 We obtain that
 $$ g_2'\px\wedge \py =-2Z_0(h)f\px\wedge \py
  = (\mp 3\lambda^2+ \mu^2 ) u\wedge v_1   + 4\lambda \mu  u\wedge v_2  \in \PH^2_\Pi.$$

Let $\ell=-1$. Then $\deg(u_j)=p, u_j=y^p , j=p+1, g_2'= 2(1-p)(\lambda x+\mu y^{p-1})xy^{p}=Z(f)$ with
$Z= (1-p)    (\lambda x+ \mu y^{p-1})y^{p-1}\px $.
  We obtain that
 $$ g_2'\px\wedge \py =\mathrm{div}(Z)f\px\wedge \py=  (1-p) \lambda  e_1f\px\wedge \py=
  \frac{1-p}{2p-1}\lambda  u\wedge v_1 \in \PH^2_\Pi.$$

Suppose that there is no $\ell\geq -1$ such that $\deg(g_2') =\rmd- \ell\cdot (\rmd-\omega_1-\omega_2)=2p-1-\ell (p-1)$.
Then   $p>2$ and $j\neq p+1, 2$.

If $j=1$, then $g_2'=(2-3p)(\lambda x+\mu y^{p-1})x=\pm (2-3p)(1-2p)\lambda u_{2p-1}+Z(f)$ for a vector field $Z$,
and we have
$$ g_2'\px\wedge \py = \pm (2-3p)(1-2p)\lambda w_{2p-1} \in \PH^2_\Pi.$$

If $3\leq j\leq p$, then $g_2'=(j-3p+1 )(\lambda x+\mu y^{p-1})xy^{j-1}= Z(f)$ for a vector field $Z$ and
$$ g_2'\px\wedge \py =0\in \PH^2_\Pi.$$

If $p+2\leq j\leq 2p-1$, then $g_2'=(j-3p+1 )(\lambda x+\mu y^{p-1})  xy^{j-1}
= Z(f)$ for a vector field $Z$ and
$$ g_2'\px\wedge \py =0\in \PH^2_\Pi.$$

If $j=2p$, $g_2'=(1-2p )(\lambda x+\mu y^{p-1})x^2
= Z(f)$ for a vector field $Z$  and
$$ g_2'\px\wedge \py =0\in \PH^2_\Pi.$$

\medskip

We have proved the following:
If $p=2$,
$$ g_2'\px\wedge \py = \left\{\begin{array}{ll} \pm 12\lambda w_{3} + 12\lambda \mu^2  u\wedge v_1+12\lambda^2 \mu\wedge v_2, &j =1, \\
    2\lambda \mu      u\wedge v_1+ 2\lambda^2  u\wedge v_2,&  j=2, \\
       - \frac{ 1}{3}\lambda  u\wedge v_1,  & j=3, \\
     (\mp 3\lambda^2+ \mu^2 ) u\wedge v_1   + 4\lambda \mu  u\wedge v_2,  & j=4;\end{array}\right.$$
if $p>2$,
$$ g_2'\px\wedge \py = \left\{\begin{array}{ll}  \pm (2-3p)(1-2p)\lambda w_{2p-1}, &j =1, \\
  \frac{6(p-1)}{2p-1}\lambda \mu      u\wedge v_1+  \frac{6(p-1)}{2p-1}\lambda^2  u\wedge v_2,& j =2, \\
      0,  & 3\leq j\leq p, \\
          \frac{1-p}{2p-1}\lambda  u\wedge v_1,  & j=p+1,\\
0& p+2\leq j\leq 2p.\end{array}\right.$$

 \medskip

We have proved the following:
If $p=2$,
$$[v_2, w_{j}]_{\sn}=\left\{\begin{array}{ll} \pm 12\lambda w_{3} -3w_{4}+ 12\lambda \mu^2  u\wedge v_1+12\lambda^2 \mu\wedge v_2, &j =1, \\
    4  \lambda \mu   u\wedge v_1+  4\lambda^2 u\wedge v_2,&  j=2, \\
        0,  & j=3, \\
    \pm 6w_{3} +   (\mp 3\lambda^2+ 5\mu^2 ) u\wedge v_1   + 8\lambda \mu  u\wedge v_2,  & j=4;\end{array}\right.$$
if $p>2$,
$$[v_2, w_{j}]_{\sn}=\left\{\begin{array}{ll}  \pm (2-3p)(1-2p)\lambda w_{2p-1}+(1-2p)w_{2p}, &j =1, \\
  \frac{12(p-1)}{2p-1} \lambda \mu u\wedge v_1+\frac{12(p-1)}{2p-1} \lambda^2 u\wedge v_2,& j =2, \\
      0,  & 3\leq j\leq 2p-1, \\
\pm p(2p-1) w_{2p-1},& j=2p.\end{array}\right.$$

\subsubsection{The final result of type $D^{\pm}_{2p}$} \label{final result for Type D2p}

We have proved the following result:
\begin{Prop} Let $\Pi=(x^2y\pm y^{2p-1})(1+\lambda x+\mu y^{p-1})\frac{\partial}{\partial x}\wedge \frac{\partial}{\partial y}$ with $p\geq 2$ be of type $D_{2p}^\pm$.
\begin{itemize}
\item[(a)]  If $p=2$,
there exists an isomorphism of Gerstenhaber  algebras
 $$\begin{array}{rll}\PH^{*}_{\Pi}&\cong &\k\langle u', v_1,v_{2}\rangle/({u'}^{2}, v_1^2,   v_{2}^{2},u' v_{1}+v_{1}u', u' v_{2}+v_{2}u', v_{1}v_{2})\\
 &&\times_{\k}\k[w_{1} ]/({w_{1}}^{2})\times_{\k}\k[w_{2}']/(w_{2}'^{2})\times_{\k}\k[w_{3}']/(w_{3}'^{2})\times_{\k}\k[w_{4}']/(w_{4}'^{2}).
 \end{array}$$
 with $|u'|=|v_1|=1, |w_{1}|=|w_{2}'|=|w_{3}'|=|w_{4}'|=2$  and the only nontrivial brackets are given by
  $$[v_1, w_1']_{\sn}=w_{2}', [v_1, w_2']_{\sn}=w_{3}', [v_1, w_4']_{\sn}=[v_2, w_2']_{\sn}=\lambda(\mu u'v_1+\lambda u'v_2), $$ $$[v_2, w_1']_{\sn}=w_{4}', [v_2, w_4']_{\sn}= 3w_{3}'+2\mu(\mu u'v_1+\lambda u'v_2).$$
 \item[(b)] If $p\geq3$,
there exists an isomorphism of Gerstenhaber  algebras
  $$\begin{array}{rll} \PH^{*}_{\Pi}&\cong& \k\langle u', v_1,v_{2}\rangle/({u'}^{2}, v_1^2,   v_{2}^{2},u' v_{1}+v_{1}u', u' v_{2}+v_{2}u', v_{1}v_{2})\\ &&\times_{\k}\k[w_{1}]/( w_{1}^{2})\times_{\k}\cdots\times_{\k}\k[w_{p-1}]/(w_{p-1}^{2})\times_{\k}\k[w_{p}']/(w_{p}'^{2})\times_{\k}\cdots\times_{\k}\k[w_{2p}']/(w_{2p}'^{2}).
  \end{array}
 $$
 with $|u'|=|v_1|=1, |w_{1}|=\cdots=|w_{p-1}|=|w_{p}'|=\cdots=|w_{2p}'|=2$ and the nontrivial brackets are given by
  $$[v_1, w_j]_{\sn}=w_{p+j-1}', 1\leq j\leq p-1, [v_1, w_p']_{\sn}=w_{2p-1}', [v_2, w_{1}]_{\sn}=w_{2p}', $$ $$[v_2, w_2]=\mu u'v_1+\lambda u'v_2, [v_2, w_{2p}']_{\sn}=\pm (1-2p)w_{2p-1}'.$$
\end{itemize}

 \end{Prop}

\begin{Proof}
(a)  Let $$\begin{array}{rcl}
w_{2}'&=& -3w_2-4\mu w_3 +4 \lambda^2\mu u\wedge v_1+4 \lambda^3 u\wedge v_2,\\
  w_{3}'&=& -3[v_1,
   w_{2}]_{\sn}) =\pm( 6w_3-3(\lambda^2 \pm     \mu^2)u\wedge v_1,\\
  w_{4}'&=&\pm 12\lambda w_{3} -3w_{4}+ 12\lambda \mu^2  u\wedge v_1+12\lambda^2 \mu\wedge v_2,\end{array}$$
 and we need to put moreover
 $u'= - 12u$.

(b)   Let $$\begin{array}{lcl}w_{p}'&=& (1-2p)w_p+ (2-3p)\mu w_{2p-1},\\
  w_{p+1}'&=& 2(1-p)w_{p+1}+ \frac{3(p-1)}{2p-1}(\lambda^2 \pm   \mu^2)u\wedge v_1,\\
  w_{p+i}'&=& (i+1-2p)w_{p+i}, 2\leq i\leq p-2,\\
  w_{2p-1}'&=&
  p(2p-1) w_{2p-1},\\
  w_{2p}'&=&  \pm (2-3p)(1-2p)\lambda w_{2p-1}+(1-2p)w_{2p},\end{array}$$ and we need to put moreover
 $u'= \frac{  12(p-1)}{2p-1} \lambda  u$.

\end{Proof}



\subsection{Type  $D_{2p+1}$:~$(x^{2}y+y^{2p})(1+\lambda x),~p\ge 2$}


In this case, $\Pi=f(1+h)\px \wedge \py$ with
$f= x^{2}y+y^{2p}, h=\lambda x,   \omega_{1}=2p-1, \omega_{2}=2, d=4p, \rmd-\omega_1-\omega_2=2p-1.$  The space $\mathcal{P}_{ \rmd-\omega_1-\omega_2}=\mathcal{P}_{ 2p-1}$ is one-dimensional with $e_1=x$ as a basis.
Since $I_{f}=\langle\frac{\partial f}{\partial x},\frac{\partial f}{\partial y}\rangle=\langle 2xy,x^{2}+2py^{2p-1}\rangle$, we can take $\{u_1=1,  \cdots, u_{2p}=y^{2p-1},  u_{2p+1}=x,  \}$ as a monomial basis of $\mathcal{M}_f=\mathcal{F}(\mathbb{K}^2)/I_{f}$, whose degrees are respectively
$0, 2, \cdots, 2(2p-1), 2p-1$.
Let $$u=(1+h)H_{f}=(1+\lambda x)((x^{2}+2py^{2p-1})\frac{\partial}{\partial x}-2xy\frac{\partial}{\partial y}), v_1=(1+h)e_1W =(1+\lambda x)x((2p-1)x\frac{\partial}{\partial x}+2y\frac{\partial}{\partial y}),$$ which is  a basis of $\PH^{1}_\Pi$, and  $w_i=u_i
\frac{\partial}{\partial x}\wedge\frac{\partial}{\partial y}$ for $1\leq i\leq 2p+1$ which form a basis of $\mathcal{M}_f\frac{\partial}{\partial x}\wedge\frac{\partial}{\partial y}$. We have an isomorphism of graded algebras  $\PH^{*}_{\Pi}\cong \wedge_\k(u, v_{1})\times_{\k}\k[w_{1}]/(w_{1}^{2})\times_{\k}\cdots\times_{\k}\k[w_{2p+1}]/(w_{2p+1}^{2}).$
  Notice that by Proposition~\ref{Prop: Wedge product}, $u\wedge v_1=4p e_1f\px\wedge \py$.

By Subsection~\ref{subsection bracket for Pi}, to compute the bracket, we need to consider
  $[v_1, w_j]_{\sn},   j=1, \cdots, 2p+1$.

  By Lemma~\ref{Lem: [v_i, w_j]},  for $1\leq j\leq 2p-1$,
$\begin{array}{rcl}
[v_1, w_j]_{\sn}&=& ((\deg(u_j)-4p)x u_j+( \deg(u_j)-6p+1) \lambda x^2 u_j )\frac{\partial}{\partial x}\wedge \frac{\partial}{\partial y}
\end{array}$
Let $g_1=(\deg(u_j)-4p)x u_j$ and $g_2=( \deg(u_j)-6p+1) \lambda x^2 u_j$.
It is easy to see that  $\deg(g_1)=2p-1+\deg(u_j)=\rmd- k\cdot (\rmd-\omega_1-\omega_2)=4p-k(2p-1)$ with $k\geq -1$ if and only if
 $k=0$ or $ 1$.

\medskip

 Let $k=0$. Then $\deg(u_j)=2p+1$, but there is no $u_j$ such that $\deg(u_j)=2p+1$. So this case is also impossible.

Let $k=1$. Then $\deg(u_j)=2, u_j=y, j=2, g_1=(2-4p)xy.$  We see that  $g_1=Y(f)$ with
$Y= (1-2p) \px$.
 Let $Y_0= (1-2p) \px,$
 and $Y_1=3\lambda (1-2p) x\px - 4\lambda y\py.$
 We obtain that
 $ g_1\px\wedge \py =-2Y_1(h)f\px\wedge \py=  \frac{  3(2p-1) }{2p}\lambda^2  u\wedge v_1\in \PH^2_\Pi.$

  Suppose that there is no $k\geq -1$ such that $\deg(g_1)=\rmd- k\cdot (\rmd-\omega_1-\omega_2)=2p-1-k(p-1)$.
Then     $j\neq 2$.
 If $j=1$, $g_1=-4px=-4p u_{2p+1}$, so $ g_1\px\wedge \py =-4p w_{2p+1}\in \PH^2_\Pi;$
 if $3\leq j\leq 2p$, $g_1=(2j-2-4p)xy^{j-1}\in I_f$, so $ g_1\px\wedge \py = 0\in \PH^2_\Pi;$
if $ j=2p+1$, then $g_1=-(1+2p )x^2=2p(1+2p)u_{2p}\, \mathrm{mod}\, I_f  $, so
 $ g_1\px\wedge \py =2p(1+2p)w_{2p}\in \PH^2_\Pi.$

\medskip

It is easy to see that  $\deg(g_2)=4p-2+\deg(u_j)=\rmd- \ell \cdot (\rmd-\omega_1-\omega_2)=4p-\ell(2p-1)$ with $\ell\geq -1$ if and only if
 $\ell=0 $ or $ -1$.

 Let $\ell=-1$. Then $\deg(u_j)=2p+1$, but there is no $u_j$ such that $\deg(u_j)=2p+1$. So this case is impossible.

Let $\ell=0$. Then $\deg(u_j)=2, u_j=y, j=2, g_2=3(1-2p)\lambda x^2y.$  We see that  $g_2=Z(f)$ with
$Z= \frac{3}{2}\lambda (1-2p) x\px$.
 Let $Z_0=-3\lambda (( 2p-1) x\px+ y\py),$
  and we obtain that
 $ g_2\px\wedge \py =-2Z_0(h)f\px\wedge \py=  \frac{3(2p-1)}{2p} \lambda^2 u\wedge v_1\in \PH^2_\Pi.$

 Suppose that there is no $\ell\geq -1$ such that $\deg(g_2)=\rmd- \ell\cdot (\rmd-\omega_1-\omega_2)=4p-\ell(2p-1)$.
Then     $j\neq 2$.
 If $j=1$, $g_2=(1-6p)\lambda x^2=2p(6p-1)\lambda u_{2p}\, \mathrm{mod}\, I_f  $, so
 $ g_2\px\wedge \py =2\lambda p(6p-1) w_{2p}\in \PH^2_\Pi;$
 if $3\leq j\leq 2p$, $g_2=(2j-6p-1)xy^{j-1}\in I_f$, so $ g_2\px\wedge \py = 0\in \PH^2_\Pi;$
if $ j=2p+1$, then $g_1=-4p\lambda x^3\in I_f  $, so
 $ g_2\px\wedge \py =0\in \PH^2_\Pi.$

\medskip

 We obtain
 $$[v_1, w_{j}]_{\sn}= \left\{\begin{array}{ll}2\lambda p(6p-1) w_{2p}-4p w_{2p+1}  , &j =1, \\
    \frac{  3(2p-1)}{p} \lambda^2  u\wedge v_1,&  j=2, \\
   0,  &  3\leq j\leq 2p, \\
     2p(1+2p)w_{2p},  & j=2p+1.\end{array}\right.$$

\medskip

We have proved the following result:
\begin{Prop} For type $D_{2p+1}$, when $\lambda=0$, then $\Pi=\Pi_0$ and
there exists an isomorphism of Gerstenhaber  algebras
     $$\PH^{*}_{\Pi_0}=\wedge_\k(u, v_{1})\times_{\k}\k[w_{1}]/(w_{1}^{2})\times_{\k}\cdots\times_{\k}\k[w_{2p-1}]/({w_{2p-1}}^{2})\times_{\k}\k[w_{2p}']/({w_{2p}'}^{2})\times_{\k}\k[w_{2p+1}']/({w_{2p+1}'}^{2}) $$
  with $|u|=|v_1|=1, |w_{1}|=\cdots=|w_{2p-1}|=|w_{2p}'|=|w_{2p+1}'|=2$  and
  the nontrivial brackets  on generators are given by
  $$[v_1, w_1]_{\sn}=w_{2p+1}', [v_1, w_{2p+1}']_{\sn}=w_{2p}'.$$

  When $\lambda\neq 0$, there exists an isomorphism of Gerstenhaber  algebras
     $$\PH^{*}_{\Pi}=\wedge_\k(u', v_{1})\times_{\k}\k[w_{1}]/(w_{1}^{2})\times_{\k}\cdots\times_{\k}\k[w_{2p-1}]/({w_{2p-1}}^{2})\times_{\k}\k[w_{2p}']/({w_{2p}'}^{2})\times_{\k}\k[w_{2p+1}']/({w_{2p+1}'}^{2}) $$
  with $|u|=|v_1|=1, |w_{1}|=\cdots=|w_{2p-1}|=|w_{2p}'|=|w_{2p+1}'|=2$  and the nontrivial brackets  on generators are given by
  $$[v_1, w_1]_{\sn}=w_{2p+1}', [v_1, w_{2p+1}']_{\sn}=w_{2p}', [v_1, w_2]=u'\wedge v_1.$$

 \end{Prop}

\begin{Proof} When $\lambda=0$,  let $w_{2p+1}'= -4p w_{2p+1}$
and $w_{2p}'= 2p(1+2p)w_{2p}$, then  $[v_1, w_1]_{\sn}=w_{2p+1}', [v_1, w_{2p+1}']_{\sn}=w_{2p}'$.
 When $\lambda\neq 0$, let $w_{2p+1}'= 2\lambda p(6p-1) w_{2p}-4p w_{2p+1}$,
  $w_{2p}'= 2p(1+2p)w_{2p}$, and we need to put moreover
 $u'= \frac{  3(2p-1)}{p} \lambda^2  u$, then $[v_1, w_2]=u'\wedge v_1$.
\end{Proof}

\subsection{Type $E_{6}$:~$x^{3}+y^{4}$}

In this type $E_6$,  $\Pi=\Pi_0=f\px\wedge \py$ with $f=x^{3}+y^{4},
   \omega_{1}=4, \omega_{2}=3, d=12, \rmd-\omega_1-\omega_2=5.$
   We    have $\mathcal{P}_{\rmd-\omega_1-\omega_2}=\mathcal{P}_5=0$ and  $\mathcal{M}_f\cong \k[x, y]/(x^2, y^{3})$, so $c=6$ and choose
   a basis of $\mathcal{M}_f $ as  $ \{u_1=1, u_2=y, u_3=x, u_4=y^2, u_5=xy,  u_6=xy^2\},$
whose degrees are respectively $0, 3, 4, 6, 7, 10$.
Let $u= H_f=4y^3\px-3x^2\py   $,
     and  $w_i=u_i \frac{\partial}{\partial x}\wedge\frac{\partial}{\partial y} $ with $1\leq i\leq 6 $, which form a basis of $\mathcal{M}_f \frac{\partial }{\partial x}\wedge \frac{\partial }{\partial y}$. 

By Theorem~\ref{Thm: wedge product result} and Subsection~\ref{subsection: bracker for Pi_0},  we obtain the following result:
\begin{Prop} \label{Prop:Type E6} Let  $\Pi_0= (x^3+y^{4})\px\wedge \py $ be Type $E_6$. We have an isomorphism of Gerstenhaber algebras  $$\PH^{*}_{\Pi_0}\cong\k[u]/(u^2)\times_{\k}\k[w_{1}]/(w_{1}^{2})\times_{\k}\cdots\times_{\k}\k[w_{6}]/(w_{6}^{2}),$$
with $|u|=1, |w_{1}|= \cdots=|w_{6 }|=2$, where the brackets are all zero.
\end{Prop}

\subsection{Type $E_{7}$:~$(x^{3}+xy^{3})(1+\lambda y^{2})$}

 In this type $E_7$, $f=x^{3}+xy^{3}, h=\lambda y^{2},
   \omega_{1}=3, \omega_{2}=2, d=9, \rmd-\omega_1-\omega_2=4.$
 We have   $\mathcal{P}_{\rmd-\omega_1-\omega_2}=\mathcal{P}_{4}=\k y^{2}$, and $r=1, e_1=y^{2}$.
 We also  have  $\mathcal{M}_f\cong \k[x, y]/(3x^2+y^3, xy^{2})$, so $c=7$ and
 choose  a monomial basis of $\mathcal{M}_f$:  $u_1=1, u_2=y, u_3=x, u_4=y^2, u_5=xy, u_6=y^3, u_7=y^4$ whose degrees are respectively
 $0, 2, 3, 4, 5, 6, 8$.

As in Theorem~\ref{Thm: wedge product result}, denote  $u=(1+h)H_f=(1+\lambda y^{2})( 3xy^2\frac{\partial}{\partial x}-(3x^2+y^3)\frac{\partial}{\partial y})$, and
    $v_1=e_1(1+h)W=y^{2}(1+\lambda y^{2})(3x\frac{\partial}{\partial x}+2y\frac{\partial}{\partial y})$ which spans   $\mathcal{P}_{\rmd-\omega_1-\omega_2}(1+h)W$,  and  $w_i=u_i \frac{\partial}{\partial x}\wedge\frac{\partial}{\partial y}$ with $1\leq i\leq 7$, which form a basis of $\mathcal{M}_f \frac{\partial }{\partial x}\wedge \frac{\partial }{\partial y}$. We have an isomorphism of graded algebras  $$\PH^{*}_{\Pi}\cong\wedge_\k(u, v_{1})\times_{\k}\k[w_{1}]/(w_{1}^{2})\times_{\k}\cdots\times_{\k}\k[w_{7}]/(w_{7}^{2}).$$

By Lemma~\ref{Lem: [v_i, w_j]},  for $1\leq j\leq 7$,
$\begin{array}{rcl}
[v_1, w_j]_{\sn}= (( \mathrm{deg}(u_{j})-9)y^{2}u_j+(\mathrm{deg}(u_{j})-13)\lambda y^4 u_j )\frac{\partial}{\partial x}\wedge \frac{\partial}{\partial y}
\end{array}$
As in Lemma~\ref{Lem: [v_i, w_j]}, let $g_1=( \mathrm{deg}(u_{j})-9)y^{2}u_j, g_2=(\mathrm{deg}(u_{j})-13)\lambda y^{4} u_j$.
Then $\deg(g_1)=\mathrm{deg}(u_{j})+4 $, $\deg(g_2)=\mathrm{deg}(u_{j})+8$.

We consider $k\in \mathbb{Z}$ such that   $\mathrm{deg}(g_1)= \rmd- k(\rmd- \omega_{1}-\omega_{2})=9-4k,$ and in this case,
$\deg(u_j)=9-4(k+1)$ and $\deg(g_2)=9-4(k+1).$
The restrictions are $  \deg(u_j)\in \{0, 2, 3, 4, 5, 6, 8\} $,  and the degrees of $g_1, g_2$ are all nonnegative.
It is easy to see that  $k=0$ and in this  case,  $\mathrm{deg}(u_j)=5, j=5, u_j=xy$,    $g_1=-4xy^3, g_2=-8 \lambda xy^5$.

We see that $g_1=Y(f)\in I_f$ with   $Y=-\frac{4}{3}y\py$ and $\mathrm{div}(Y)=  -\frac{4}{3}$. Since $\mathrm{deg}(g_1)=9=\rmd- 0\cdot (\rmd- \omega_{1}-\omega_{2})$,   by Proposition~\ref{Prop: final computation},
 let $Y_0= -x\px-2y\py.$
 Then $$[g_1\px\wedge \py, w_{5}]_{\sn}=-2Y_0(h)f\px\wedge \py= 8\lambda e_1f\px\wedge \py=\frac{8}{9}\lambda u\wedge v_1\in \PH^2_{\Pi}. $$
 Since $\mathrm{deg}(g_2)=13=\rmd- (-1)\cdot (\rmd- \omega_{1}-\omega_{2})$, and
  $g_2=Z(f)\in I_f$ with $ Z=-\frac{8}{3}\lambda y^3\py$ and $ \mathrm{div}(Z)= -8 \lambda y^{2}$, by Proposition~\ref{Prop: final computation}, $$[g_2\px\wedge \py, w_{5}]_{\sn}= \mathrm{div}(Z)f\px\wedge \py=  -8 \lambda e_1f\px\wedge \py= -\frac{8}{9}\lambda  u\wedge v_1\in \PH^2_{\Pi}. $$
We obtain that $[v_1, w_5]_{\sn}=0\in \PH^2_{\Pi}.$

Now for $j\neq 5$, using Proposition~\ref{Prop: final computation},  we obtain easily that
$$\begin{array}{rclcl}
  [v_1, w_1]_{\sn}  &=&  (-9y^2-13\lambda y^4)\px\wedge \py  &=&-9w_4-13\lambda w_7\in \PH^2_\Pi,\\

    [v_1, w_2]_{\sn} &=&  (-7y^3-11\lambda y^5)\px\wedge \py  &=& -7w_6\in \PH^2_\Pi, \\

   [v_1, w_3]_{\sn}  &=&  (-6xy^2-10\lambda xy^4)\px\wedge \py&=& 0\in \PH^2_\Pi,\\

     [v_1, w_4]_{\sn} &=&  (-5y^4-9\lambda  y^6)\px\wedge \py &=& -5w_7\in \PH^2_\Pi, \\

    [v_1, w_6]_{\sn} &=&   (-3y^5-7\lambda  y^7)\px\wedge \py &=& 0\in \PH^2_\Pi, \\

   [v_1, w_7]_{\sn}  &=&   (- y^6-5\lambda  y^8)\px\wedge \py &=& 0\in \PH^2_\Pi. \end{array}$$

\begin{Prop} \label{Prop: Type E7} For type $E_7$, as
  Gerstenhaber  algebras,  $\PH^{*}_{\Pi}$  is isomorphic to
     $$\wedge_\k(u, v_{1})\times_{\k}\k[w_{1}]/(w_{1}^{2})\times_{\k}\k[w_{2}]/(w_{2}^{2})\times_{\k}\k[w_{3}]/(w_{3}^{2})
     \times_{\k}\k[w_{4}']/({w_{4}'}^{2})\times_{\k}\k[w_{5}]/(w_{5}^{2})\times_{\k}\k[w_{6}']/({w_{6}'}^{2})
     \times_{\k}\k[w_{7}']/({w_{7}'}^{2}) $$
  with $|u|=|v_1|=1, |w_{1}|=|w_{2}|=|w_{3}|=|w_{4}'|=|w_{5}|=|w_{6}'|= |w_7'|=2$  and   the nontrivial brackets on generators are given by
  $$[v_1, w_1]_{\sn}=w_4', [v_1, w_2]_\sn=w_6', [v_1, w_4']_\sn=w_7'.$$

 \end{Prop}

 \begin{Proof} It suffices to define $w_4'= -9w_4-13\lambda w_7$, $w_6'=-7w_6$ and $w_7'=45 w_7$.

 \end{Proof}

 \subsection{Type $E_{8}$:~$x^{3}+y^{5}$}

 In this type $E_8$, $\Pi=\Pi_0=f\px\wedge \py, $ with $ f=x^{3}+y^{5},
   \omega_{1}=5, \omega_{2}=3, d=15, \rmd-\omega_1-\omega_2=7.$
   We    have $\mathcal{P}_{\rmd-\omega_1-\omega_2}=\mathcal{P}_7=0$ and  $\mathcal{M}_f\cong \k[x, y]/(x^2, y^{4})$, so $c=8$ and choose
   a basis of $\mathcal{M}_f $ as  $ u_1=1, u_2=y, u_3=x, u_4=y^2, u_5=xy, u_6=y^3, u_7=xy^2, u_8=xy^3.$
Let $u= H_f=5y^4\px-3x^2\py   $,
     and  $w_i=u_i \frac{\partial}{\partial x}\wedge\frac{\partial}{\partial y}$ with $1\leq i\leq 8 $, which form a basis of $\mathcal{M}_f \frac{\partial }{\partial x}\wedge \frac{\partial }{\partial y}$.
 We have  the following result:
\begin{Prop} \label{Prop:Type E8} Let  $\Pi_0= (x^3+y^{5})\px\wedge \py $ be Type $E_8$. We have an isomorphism of Gerstenhaber algebras  $$\PH^{*}_{\Pi_0}\cong \k[u]/(u^2)\times_{\k}\k[w_{1}]/(w_{1}^{2})\times_{\k}\cdots\times_{\k}\k[w_{8}]/(w_{8}^{2}),$$
where the brackets are all zero.
\end{Prop}






%



\subsection{Final result for simple singularities}\label{Section: Final result for simple singularities}
We summarize the computation for simple singularities in Table~\ref{Simple singularity computations}.
\begin{table}[h]
 \centering
\begin{tabular}{|c|c|}
\hline
 type  & $\mathrm{HP}^{*}$ \\
\hline
$A_{2p}$ & $\wedge_\k(u, v_{1})\times_{\mathbb{K}}\mathbb{K}[w_{1}]/(w_{1}^{2})
\times_{\mathbb{K}}\cdots\times_{\mathbb{K}}\mathbb{K}[w_{2p}]/(w_{2p}^{2}) $\\
$p\ge 1$&$|u|=|v_1|=1,  |w_{1}|=|w_{2}|=\cdots=|w_{2p}|=2$  \\

\hline
$A_{1}^{\pm}$ & $\wedge_\k(u, v_1')\times_{\mathbb{K}}
 \mathbb{K}[w_1]/(w_{1}^{2}) $ \\
   &$|u|=|v_1'|=1, |w_{1}|=2$  \\
 & $[v_1', w_{1}]_{\sn}=  w_1 $\\
  \hline
$A_{3}^{\pm}$ & $\wedge_\k(u, v_{1})\times_{\k}\k[w_{1}]/(w_{1}^{2})\times_{\k}\k[w_{2}']/{(w_{2}'}^{2}) \times_{\k}\k[w_{3}']/({w_{3}'}^{2}) $ \\
$\lambda=0 $  &$|u|=|v_1|=1, |w_{1}|= |w_2'|= |w_{3}'|=2$  \\
&$[v_1, w_1]_{\sn}=w_2', [v_1, w_2']_{\sn}=w_3' $\\

  \hline
  $A_{3}^{\pm}$ & $ \wedge_\k(u', v_{1})\times_{\k}\k[w_{1}]/(w_{1}^{2}) \times_{\k}\k[w_{2}']/{(w_{2}'}^{2}) \times_{\k}\k[w_{3}']/({w_{3}'}^{2})
  $ \\
$\lambda\neq 0 $  &$ |u'|=|v_1|=1, |w_{1}|=|w_2'|=|w_3'|=2$  \\
&$ [v_1, w_1]_{\sn}=w_2', [v_1, w_2']_{\sn}=w_3', [v_1, w_3']=u'v_1 $\\

  \hline
$A_{2p-1}^{\pm}$ & $\wedge_\k(u, v_{1})\times_{\k}\k[w_{1}]/(w_{1}^{2})\times_{\k} \cdots \times_{\k}\k[w_{p-1}]/(w_{p-1}^{2})
     \times_{\k}$\\
  $p\ge 3$    & $\k[w_{p}']/({w_{p}'}^{2})\times_{\k}  \cdots \times_{\k}\k[w_{2p-1}']/({w_{2p-1}'}^{2}) $ \\
 & $|u|=|v_1|=1, |w_{1}|=\cdots=|w_{p-1}|=|w_p'|=\cdots=|w_{2p-1}'|=2$, \\
&$ [v_1, w_1]_\sn=w_p', [v_1, w_p']_\sn=w_{2p-1}', [v_1, w_j]_{\sn}=w_{p+j-1}',j=2, \cdots, p-1 $\\
\hline
$D_{4}^{\pm}$ &
  $\k\langle u', v_1,v_{2}\rangle/({u'}^{2}, v_1^2,   v_{2}^{2},u' v_{1}+v_{1}u', u' v_{2}+v_{2}u', v_{1}v_{2})\times_{\k} $ \\
 & $\k[w_{1} ]/(w_{1}^{2})\times_{\k}\k[w_{2}']/(w_{2}'^{2})\times_{\k}\k[w_{3}']/(w_{3}'^{2})\times_{\k}\k[w_{4}']/(w_{4}'^{2}) $  \\
  & $|u'|=|v_1|=|v_{2}|=1, |w_{1}|= |w_2'|=|w_3'|=|w_{4}|=2$   \\
& $[v_1, w_1']_{\sn}=w_{2}', [v_1, w_2']_{\sn}=w_{3}', [v_1, w_4']_{\sn}=[v_2, w_2']_{\sn}=\lambda(\mu u'v_1+\lambda u'v_2),$\\
& $ [v_2, w_1']_{\sn}=w_{4}', [v_2, w_4']_{\sn}= \mp 3w_{3}'+2\mu(\mu u'v_1+\lambda u'v_2)$ \\

\hline
$D_{2p}^{\pm}$ &
  $\k\langle u', v_1,v_{2}\rangle/({u'}^{2}, v_1^2,   v_{2}^{2},u' v_{1}+v_{1}u', u' v_{2}+v_{2}u', v_{1}v_{2})\times_{\k}$\\
 $p\ge 3$   & $\k[w_{1}]/( w_{1}^{2})\times_{\k}\cdots\times_{\k}
  \k[w_{p-1}]/(w_{p-1}^{2})\times_{\k}\k[w_{p}']/(w_{p}'^{2})\times_{\k}\cdots\times_{\k}\k[w_{2p}']/({w_{2p}'}^{2})  $  \\
 & $|u'|=|v_1|=1, |w_{1}|=\cdots=|w_{p-1}|=|w_{p}'|=\cdots=|w_{2p}'|=2$   \\
& $[v_1, w_j]_{\sn}=w_{p+j-1}', 1\leq j\leq p-1, [v_1, w_p']_{\sn}=w_{2p-1}', [v_2, w_{1}]_{\sn}=w_{2p}',$ \\
& $[v_2, w_2]=\mu u'v_1+\lambda u'v_2, [v_2, w_{2p}']_{\sn}=\pm (1-2p)w_{2p-1}'$\\

\hline $D_{2p+1}$ &  $\wedge_\k(u, v_{1})\times_{\k}\k[w_{1}]/(w_{1}^{2})\times_{\k}\cdots\times_{\k}\k[w_{2p-1}]/({w_{2p-1}}^{2})
\times_{\k}$ \\
$p\ge 2$  & $\k[w_{2p}']/({w_{2p}'}^{2})\times_{\k}\k[w_{2p+1}']/({w_{2p+1}'}^{2}) $ \\
 $\lambda=0$& $ |u|=|v_1|=1, |w_{1}|=\cdots=|w_{2p-1}|=|w_{2p}'|=|w_{2p+1}'|=2 $  \\
  &$[v_1, w_1]_{\sn}=w_{2p+1}', [v_1, w_{2p+1}']_{\sn}=w_{2p}'$ \\

\hline $D_{2p+1}$& $\wedge_\k(u', v_{1})\times_{\k}\k[w_{1}]/(w_{1}^{2})\times_{\k}\cdots\times_{\k}\k[w_{2p-1}]/({w_{2p-1}}^{2})
\times_{\k}$\\
$p\ge 2$ &$\k[w_{2p}']/({w_{2p}'}^{2})\times_{\k}\k[w_{2p+1}']/({w_{2p+1}'}^{2}) $
 \\ $\lambda\neq 0 $ &  $|u|=|v_1|=1, |w_{1}|=\cdots=|w_{2p-1}|=|w_{2p}'|=|w_{2p+1}'|=2$    \\
 & $[v_1, w_1]_{\sn}=w_{2p+1}', [v_1, w_{2p+1}']_{\sn}=w_{2p}', [v_1, w_2]=u'  v_1$   \\

\hline

$E_{6}$& $\mathbb{K}[u]/(u^{2})\times_{\mathbb{K}}\mathbb{K}[w_{1}]/(w_{1}^{2})
\times_{\mathbb{K}}\cdots\times_{\mathbb{K}}\mathbb{K}[w_{6}]/(w_{6}^{2}),$\\ &$|u|=1,~|w_{1}|= \cdots=|w_{6}|=2$   \\

\hline
$E_{7}$ & $\wedge_\k(u, v_{1})\times_{\k}\k[w_{1}]/(w_{1}^{2})\times_{\k} \cdots
     \times_{\k}\k[w_{3}]/(w_{31}^{2})\times_{\k}\k[w_{4}']/({w_{4}'}^{2})\times_{\k}$\\
     &$\k[w_{5}]/(w_{5}^{2})\times_{\k}\k[w_{6}']/({w_{6}'}^{2})
     \times_{\k}\k[w_{7}']/({w_{7}'}^{2}) $\\
 &$|u|=|v_1|=1, |w_{1}|=|w_{2}|=|w_{3}|=|w_{4}'|=|w_{5}|=|w_{6}'|= |w_7'|=2 $   \\
 &$[v_1, w_1]_{\sn}=w_4', [v_1, w_2]_\sn=w_6', [v_1, w_4']_\sn=w_7'   $\\

\hline
$E_{8}$ & $\mathbb{K}[u]/(u^{2})\times_{\mathbb{K}}
\mathbb{K}[w_{1}]/(w_{1}^{2})\times_{\mathbb{K}}\cdots\times_{\mathbb{K}}\mathbb{K}[w_{8}]/(w_{8}^{2}) $\\ &$|u|=1, |w_{1}|= \cdots=|w_{8}|=2$  \\

\hline
\end{tabular}
\caption{Gerstenaber algebra structure for simple singularities}
\label{Simple singularity computations}
 \end{table}

\medskip

\textbf{Acknowledgements}  The  authors were supported by Natural Science Foundation of China (No. 11671139, 11971460)   and by Science and Technology Commission of Shanghai Municipality (No. 18dz2271000).

We are very grateful to the comments of the referee(s) and the editor which led to a substantial revision of the text. 

\textbf{Data availability} Data sharing is not applicable to this article as no new data were created or analyzed in this study.

\end{document}